%% file: chapter.tex
\documentclass[11pt]{book}

\usepackage[utf8]{inputenc}
\usepackage[T1]{fontenc}
\usepackage{lmodern}

\usepackage{amsmath,amssymb,amsthm,mathtools}

\usepackage{booktabs}
\usepackage{array}
\usepackage{multirow}

\usepackage{graphicx}
\usepackage{tikz}
\usetikzlibrary{arrows.meta,fit,positioning,decorations.pathreplacing}
\usepackage{caption}
\usepackage{hyperref}
\usepackage[capitalise,noabbrev]{cleveref}

\providecommand{\doi}[1]{doi: \href{https://doi.org/#1}{\begingroup\urlstyle{rm}\nolinkurl{#1}\endgroup}}

\usepackage[numbers,sort&compress]{natbib}

\theoremstyle{definition}
\newtheorem{definition}{Definition}[section]
\newtheorem{example}{Example}[section]

\theoremstyle{plain}
\newtheorem{theorem}{Theorem}[section]

\theoremstyle{remark}
\newtheorem{remark}{Remark}[section]

\newcommand{\R}{\mathbb{R}}
\newcommand{\Z}{\mathbb{Z}}

\newcommand{\Dgm}{\mathrm{Dgm}}
\newcommand{\dB}{d_B}
\newcommand{\kmax}{k\text{-}\max}
\newcommand{\NN}{\mathrm{NN}}

\usepackage[margin=1in]{geometry}

\begin{document}

\label{ch:tda}

\begingroup
\centering
{\LARGE\bfseries Persistent Homology of Time Series\\through Complex Networks\par}
\vspace{0.8em}
{\normalsize İsmail Güzel\par}
\vspace{0.3em}
{\small\itshape Institute of Applied Mathematics, Middle East Technical University, Ankara, Türkiye\par}
{\small\itshape Network Technologies Department, TÜBİTAK ULAKBİM, Ankara, Türkiye\par}
\vspace{0.2em}
{\small\ttfamily iguzel@metu.edu.tr\par}
\vspace{1.2em}
\endgroup

\begin{quote}\small
\textbf{Abstract.}
We present a unified pipeline for univariate time series classification via complex networks and persistent homology.  A time series is mapped to a graph through one of five constructions across three families---visibility (natural and horizontal visibility graphs), transition, and proximity---and the graph is converted to a dissimilarity matrix from which a Vietoris--Rips filtration yields persistence diagrams.  These diagrams are vectorized into fixed-length features through persistence landscapes and topological summary statistics.  By standardizing the downstream processing, differences in classification performance are attributable to the network construction and distance metric alone.  Experiments on twelve UCR benchmarks show that (i)~no single construction dominates: the optimal graph type depends on the signal's discriminative structure; (ii)~the graph distance metric is a first-order design choice, with diffusion distance uniformly outperforming shortest-path alternatives; and (iii)~persistence-based features degrade gracefully under noise, consistent with the classical stability theorem of persistent homology.

\end{quote}
\bigskip

\section{Introduction}
\label{sec:intro}

\subsection{Why Topology for Time Series?}
\label{sec:motivation}

Time series data appear across virtually every scientific and engineering domain, and reliable classification remains a central problem in modern data analysis~\cite{Bagnall2017}.  Classical approaches remain essential, but each emphasizes only part of the structure present in a signal.  Frequency-domain methods are most effective when the underlying process is sufficiently stationary or can be localized meaningfully in frequency, conditions that are often strained in strongly nonstationary applications.  Recent financial-forecasting studies likewise show that richer topological features can improve predictive performance when standard representations miss nonlinear structure~\cite{Carlos2025}.  Linear autoregressive models provide interpretable baselines, yet they are naturally best suited to dynamics that can be approximated well through linear dependence.  Dynamic time warping offers flexible pairwise alignment, but its pairwise cost grows quadratically with the series length and can become computationally expensive at scale~\cite{Bagnall2017}.  Wavelet-based methods address nonstationarity more directly and often provide competitive baselines, though recent comparative studies show that topological features can capture distinctions that standard wavelet and descriptive features miss in stochastic-process classification~\cite{Guzel2023}.  These complementary strengths and limitations motivate representations that capture the \emph{global geometric structure} of a signal rather than relying exclusively on local statistical summaries.

Topological Data Analysis (TDA) offers a natural response to this need.  Its central engine, persistent homology, tracks how topological features---connected components, loops, and voids---appear and disappear as one varies a scale parameter across the data.  The resulting persistence diagram is a compact multi-scale summary of shape that avoids committing to a single analysis scale \emph{a priori}~\cite{DeyWang2022}.  A key theoretical guarantee supports this viewpoint: the stability theorem of Cohen-Steiner, Edelsbrunner, and Harer~\cite{CohenSteiner2007} shows that small perturbations in the input induce only small changes in the persistence diagram.  In time-series applications, topological features have therefore emerged as a useful complement to more classical descriptors.  Recent studies include classification and clustering pipelines for univariate time series, including simulated models, financial data, and physiological-signal datasets~\cite{Majumdar2020,Karan2021}, tutorial treatments focused specifically on persistent homology for time series~\cite{Ravishanker2021}, persistent-homology analysis of time-series distance matrices~\cite{Ichinomiya2023}, feature-aware formulations that incorporate domain knowledge into persistence calculations~\cite{Heo2024}, machine-learning pipelines built on persistent-homology summaries of time series~\cite{Ichinomiya2025}, and forecasting frameworks that inject topological features into modern predictive models~\cite{Lin2026}.

Applying TDA to a univariate time series requires bridging a representational gap: the one-dimensional sequence must be lifted into a geometric object on which persistent homology can act.  A standard way to do so is through Takens' delay embedding~\cite{Takens1981}, which reconstructs the state space by forming vectors of lagged observations.  Although the embedding itself is inexpensive, the subsequent Vietoris--Rips computation on the resulting point cloud can become prohibitively expensive for long signals, motivating efficient implementations such as Ripser~\cite{Bauer2021}.

Complex networks provide an attractive alternative.  Mapping a time series into a graph---through visibility algorithms~\cite{Lacasa2008}, ordinal partition networks~\cite{Myers2019_PRE,Myers2023_SIAM}, coarse-grained state-space networks~\cite{Myers2023_CGSSN}, or proximity graphs such as $k$-nearest neighbor graphs---can preserve salient dynamical information in a substantially smaller combinatorial object.  The subsequent topological computation is then carried out on a graph-derived pairwise matrix rather than on the full embedded point cloud.  More generally, persistent homology has been used to study complex networks directly~\cite{Horak2009}, so graph-based time-series representations provide a natural intermediary between signal data and topological descriptors.  The resulting pipeline, \emph{time series $\to$ complex network $\to$ persistence diagram}, retains the multi-scale and noise-robust perspective of TDA while exploiting the compactness and interpretability of network representations~\cite{Silva2021}.

\subsection{Scope and Contributions}
\label{sec:scope}

We develop a unified pipeline for univariate time series classification that routes through complex networks and persistent homology.  We organize five network construction methods into three families, following the taxonomy of Silva et~al.~\cite{Silva2021}: \emph{visibility-based} methods (natural visibility graph and horizontal visibility graph, evaluated as two separate pipelines), \emph{transition-based} methods (ordinal partition networks and coarse-grained state-space networks), and \emph{proximity-based} methods ($k$-nearest neighbor graphs).  For each network type, we examine several graph-based separation measures that convert the network into a symmetric pairwise matrix suitable for the downstream flag filtration.

To enable a fair comparison across methods, we \emph{standardize} the vectorization stage: every persistence diagram is mapped to a fixed-length feature vector using persistence landscapes~\cite{Bubenik2015}, Adcock--Carlsson coordinates, and diagram statistics such as persistent entropy and total persistence.  This design isolates the effect of the network construction and graph-based separation measure from the downstream learning step.

We also delimit the scope.  We restrict attention to univariate time series; multivariate extensions are surveyed elsewhere.  We do not treat sublevel-set filtrations applied directly to the time series function~\cite{Myers2022_ANAPT}, sliding-window persistence for periodicity detection~\cite{Perea2015}, or learning-based forecasting models augmented with topological features~\cite{Lin2026}.

The remainder is organized as follows.  \Cref{sec:common-route} develops the standardized downstream route shared by all methods---from simplicial complexes and persistent homology through graph distance metrics to persistence diagram extraction.  \Cref{sec:networks} introduces the five network construction methods in detail.  \Cref{sec:vectorization} covers vectorization and featurization of persistence diagrams.  \Cref{sec:experiments} validates the pipeline empirically on twelve UCR benchmarks through three experiments: a graph-construction comparison, a distance-matrix ablation, and a noise-robustness study.  \Cref{sec:conclusion} offers conclusions, practical recommendations, and open problems.

\section{From Complex Networks to Persistence Diagrams}
\label{sec:common-route}

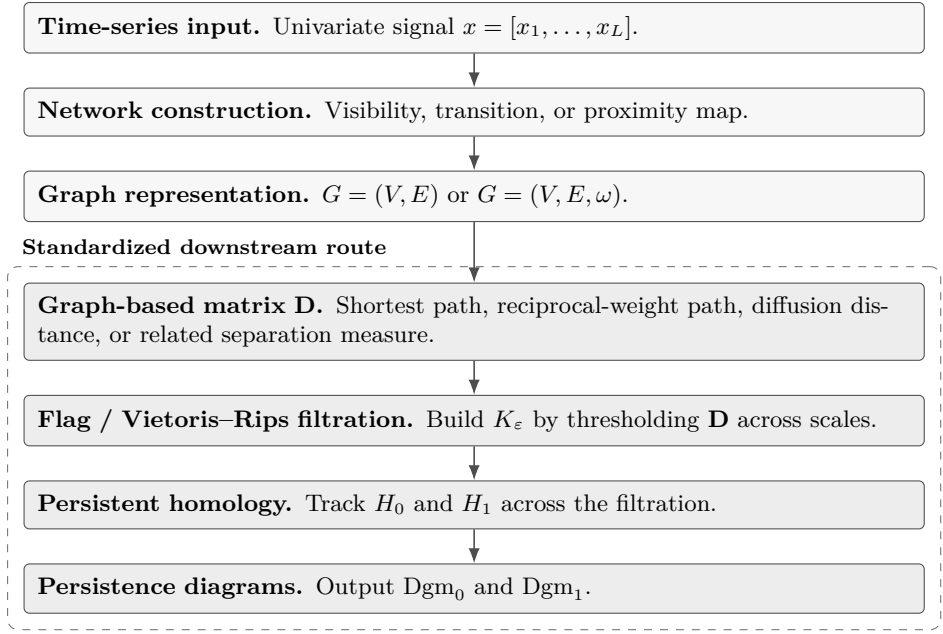
\begin{figure}[t]
\centering
\input{figures/common-route-flow}
\caption{Schematic for \cref{sec:common-route}.  Once a network has been constructed from the time series, the remaining map from $G$ to $\Dgm_0$ and $\Dgm_1$ is held fixed across methods.}
\label{fig:common-route-flow}
\end{figure}

Regardless of which network construction method is employed (\cref{sec:networks}), the downstream path from a complex network to a persistence diagram follows a single standardized route.  This section develops that route from first principles: we introduce the necessary algebraic and combinatorial tools, define the graph-based matrices used downstream, build the Vietoris--Rips filtration, and extract persistence diagrams.  Fixing these steps across all network types is what allows us to attribute differences in classification performance to the network construction itself.  The flow chart in \cref{fig:common-route-flow} summarizes this route before the subsections develop each stage in detail.  Readers seeking a comprehensive treatment of the underlying topology are referred to~\cite{DeyWang2022}.

We begin by fixing the graph input and notation.  The common route begins with the undirected graph obtained from a univariate time series by one of the methods described in \cref{sec:networks}; for transition networks, this means after symmetrization when needed.  In this study, the corresponding \emph{complex network} is modeled mathematically by that graph: we use \emph{network} when emphasizing the signal-derived representation and \emph{graph} when emphasizing the combinatorial object carried into the downstream topology.  We write $G = (V, E)$ in the unweighted case and $G = (V, E, \omega)$ when positive edge weights $\omega$ are available.  Following the taxonomy of Silva et~al.~\cite{Silva2021}, these mappings from time series to networks fall into three families---visibility, transition, and proximity---which encode different structural aspects of the signal, including visibility relations, transition frequencies, and local geometric proximity in the reconstructed state space.  With this input graph fixed, the next step is to pass from $G$ to the simplicial structures on which persistent homology is computed.

\subsection{Simplicial Complexes and Filtrations}
\label{sec:simplicial}

\begin{definition}[Simplicial complex]\label{def:simplicial-complex}
A \emph{simplicial complex}~$K$ on a finite vertex set~$V$ is a collection of subsets of~$V$, called \emph{simplices}, that is closed under taking subsets: if $\sigma \in K$ and $\sigma' \subseteq \sigma$, then $\sigma' \in K$.  A simplex with $k+1$ vertices is a \emph{$k$-simplex}: $0$-simplices are vertices, $1$-simplices are edges, $2$-simplices are filled triangles, and so on.
\end{definition}

A simplicial complex of dimension~$1$ is simply a graph; the passage to higher dimensions is what allows TDA to detect features---such as loops and cavities---that graphs alone cannot capture~\cite{Chazal2021}.  A \emph{filtration} is a nested sequence of subcomplexes $K_0 \subseteq K_1 \subseteq \cdots \subseteq K_m = K$ indexed by an increasing parameter~$\varepsilon$.  As $\varepsilon$ grows, simplices are added but never removed~\cite{Edelsbrunner2002}.

\subsection{Constructing a Graph-Based Pairwise Matrix}
\label{sec:distance-matrix}

Let $G=(V,E)$ be a finite undirected graph.  When positive edge weights are available, we write $G=(V,E,\omega)$, where $\omega(e)>0$ is the weight of edge~$e$.  For vertices $a,b\in V$, a path from $a$ to $b$ is a sequence of adjacent vertices $P=(v_0=a,v_1,\dots,v_m=b)$; its hop count is $|P|=m$.  When weights encode similarity or transition frequency, stronger edges should correspond to smaller separations, so we convert edge weights to edge lengths by setting $\ell(e)=1/\omega(e)$.

The next stage of the pipeline uses a symmetric nonnegative matrix $\mathbf{D} \in \R^{|V|\times |V|}$ with zero diagonal.  When $\mathbf{D}$ satisfies the metric axioms, it defines a finite metric space on~$V$; when it does not, it is still a graph-derived dissimilarity matrix that can be thresholded by the same pairwise rule.  The adjacency matrix alone is insufficient because it records only immediate neighbors rather than multistep connectivity.  We therefore replace $G$ by one of the four matrices below.  The first uses only the unweighted adjacency structure, whereas the latter three use edge weights when those weights carry meaningful information.

\paragraph{Shortest unweighted path.}
$\mathbf{D}(a,b)$ is the minimum hop count over all paths from $a$ to $b$:
\[
  \mathbf{D}(a,b) = \min_{P:a\leadsto b} |P|.
\]
This is the usual shortest-path distance in the unweighted graph, computed via breadth-first search.  It captures the graph's combinatorial structure but ignores edge weights.

\paragraph{Hop count on frequency-optimal path.}
Among all paths that minimize total reciprocal weight, we report the smallest hop count:
\[
  \mathbf{D}(a,b)
  = \min\Bigl\{|P| : P \in \arg\min_{Q:a\leadsto b} \sum_{e \in Q} \ell(e)\Bigr\}.
\]
This variant measures how many state transitions separate two nodes along a route that is optimal with respect to edge frequencies.  Because it mixes two optimization criteria, it is best viewed as a graph-derived \emph{dissimilarity} rather than as a guaranteed metric.

\paragraph{Reciprocal-weight shortest path.}
Here edge weights are first inverted to obtain dissimilarities, and the distance is defined as the minimum total dissimilarity:
\[
  \mathbf{D}(a,b) = \min_{P:a\leadsto b} \sum_{e \in P} \ell(e).
\]
High-frequency transitions contribute small terms, so well-traveled routes yield short distances.  This is the standard weighted shortest-path metric on the reciprocal-weight graph, computed via Dijkstra's algorithm.

\paragraph{Diffusion distance.}
Introduced by Coifman and Lafon~\cite{Coifman2006}, this measures multiscale connectivity through random walks.  Let $\mathbf{A}$ denote the adjacency matrix, binary in the unweighted case and weighted when edge weights are available, define the row-stochastic transition matrix by $\mathbf{P}(i,j)=\mathbf{A}(i,j)/\sum_k \mathbf{A}(i,k)$, and use the lazy walk $\tilde{\mathbf{P}}=\tfrac{1}{2}(\mathbf{P}+\mathbf{I})$.  Writing
\[
  s(c)=\sum_j \mathbf{A}(c,j),
  \qquad
  \pi(c)=\frac{s(c)}{\sum_k s(k)}.
\]
The diffusion distance after $t$ steps is
\begin{equation}\label{eq:diffusion}
  d_t(a,b) = \sqrt{\sum_{c \in V} \frac{1}{\pi(c)}\bigl[\tilde{\mathbf{P}}^t(a,c) - \tilde{\mathbf{P}}^t(b,c)\bigr]^2}.
\end{equation}
Two nodes connected by many high-probability paths have small diffusion distance, even if a spurious low-weight edge exists between them.  This quantity isolates long-range structural similarities by comparing random-walk distributions rather than a single optimal path.  In the weighted OPN experiments of Myers et~al.~\cite{Myers2023_SIAM}, a practical heuristic was to choose $t$ in the range $\mathrm{diam}(G) < t < 3\,\mathrm{diam}(G)$, where $\mathrm{diam}(G)$ denotes the unweighted graph diameter; we retain that guideline here as a tuning heuristic rather than as a theorem.

These four constructions form a practical spectrum.  Unweighted shortest path is the most interpretable baseline when only connectivity matters.  Hop count on the frequency-optimal path preserves an integer notion of separation while still using edge weights to choose the route.  Reciprocal-weight shortest path is often the natural default when edge weights carry substantive information and one wants a genuine path-length metric.  Diffusion distance is the most global and typically the most computationally demanding option, but it is often more robust to isolated low-weight edges because it compares full random-walk profiles rather than a single optimal path.  Normalizing $\mathbf{D}$ to $[0,1]$ further facilitates comparison across datasets and network constructions.  Once $\mathbf{D}$ has been chosen, the next subsection turns it into a filtration by thresholding its pairwise entries.

\subsection{The Vietoris--Rips Filtration}
\label{sec:rips-computation}

At this stage, the graph itself is no longer the direct input; the filtration is built from the symmetric pairwise matrix $\mathbf{D}$ constructed above.  For an unweighted graph, the default choice is the shortest-path matrix in hop count.  For a weighted graph, the default depends on the meaning of the weights: if larger weights encode stronger similarity or more frequent transitions, one first sets $\ell(e)=1/\omega(e)$ and uses reciprocal-weight shortest path; if the weights already represent costs or distances, one uses the weighted shortest path directly.  Diffusion distance remains an optional more global alternative in either case.  In particular, the adjacency matrix itself is not the object that is thresholded in the filtration.

Given the matrix $\mathbf{D}$, we threshold it at increasing scales to construct a nested sequence of simplicial complexes.  For a monotonically increasing scale parameter $\varepsilon \in \R_{\geq 0}$, we define
\begin{equation}\label{eq:rips-from-graph}
  K_\varepsilon = \bigl\{\sigma \subseteq V \mid \mathbf{D}(u,v) \leq \varepsilon \;\;\text{for all}\; u \neq v \in \sigma\bigr\}.
\end{equation}
Thus a finite vertex set $\sigma=\{v_0,\dots,v_k\}$ is included exactly when every pair of its vertices satisfies $\mathbf{D}(v_i,v_j)\leq\varepsilon$.  Because a simplex is added wherever a fully connected subgraph (clique) exists at scale~$\varepsilon$, $K_\varepsilon$ is a \emph{clique complex} (or \emph{flag complex}), entirely determined by its $1$-skeleton~\cite{Bauer2021}.  This means one needs only pairwise matrix entries to reconstruct the full complex---no ambient coordinates are required.  When $\mathbf{D}$ is a metric, $K_\varepsilon$ is exactly the \emph{Vietoris--Rips complex} at scale~$\varepsilon$; for a general symmetric graph-derived dissimilarity matrix, the same formula defines the corresponding flag filtration.  As $\varepsilon$ increases, we obtain the required filtration $K_{\varepsilon_1} \subseteq K_{\varepsilon_2} \subseteq \cdots \subseteq K_{\varepsilon_m} = K$.

\subsection{Persistent Homology}
\label{sec:homology}

To extract topological information from the filtration, we compute homology at each scale $\varepsilon$.  We work over the field $\Z_2 = \{0,1\}$, which is standard in computational TDA because arithmetic is mod~2 and one need not keep track of simplex orientations.  For a fixed complex $K_\varepsilon$, the $p$-chain group $C_p(K_\varepsilon)$ is the vector space spanned by its $p$-simplices.  The \emph{boundary operator}
\[
  \partial_p \colon C_p(K_\varepsilon) \to C_{p-1}(K_\varepsilon)
\]
maps each $p$-simplex to the sum of its $(p{-}1)$-dimensional faces.  The identity $\partial_{p-1}\circ\partial_p=0$ implies that every boundary is automatically a cycle.  Writing
\[
  Z_p(K_\varepsilon)=\ker(\partial_p),
  \qquad
  B_p(K_\varepsilon)=\operatorname{im}(\partial_{p+1}),
\]
the $p$-th homology group is the quotient
\begin{equation}\label{eq:homology}
  H_p(K_\varepsilon) = Z_p(K_\varepsilon)\,\big/\,B_p(K_\varepsilon).
\end{equation}
Its dimension
\[
  \beta_p(K_\varepsilon)=\dim H_p(K_\varepsilon)
\]
is the \emph{$p$-th Betti number} at scale~$\varepsilon$.  Intuitively, $\beta_0$ counts connected components, $\beta_1$ counts independent loops, and $\beta_2$ counts enclosed cavities.  We focus on $H_0$ and $H_1$.

Persistent homology records how these groups change as $\varepsilon$ increases.  Whenever $\varepsilon_i \leq \varepsilon_j$, the inclusion $K_{\varepsilon_i}\subseteq K_{\varepsilon_j}$ induces a linear map
\[
  H_p(K_{\varepsilon_i}) \to H_p(K_{\varepsilon_j}).
\]
A homology class is \emph{born} at parameter $b$ when it first appears in the filtration.  It \emph{dies} at parameter $d$ when it ceases to represent an independent class at a later scale.  For $H_0$, this happens when two connected components merge and the younger class is paired away; for $H_1$, it happens when enough higher-dimensional simplices enter to make the loop a boundary.  The quantity $d-b$ is the \emph{persistence} of the feature.

The collection of lifespans is summarized by the \emph{persistence diagram} $\Dgm_p$, a multiset of birth-death pairs $\{(b_i,d_i)\}_{i\in I}$.  Features far from the diagonal persist across many scales and therefore indicate stronger multi-scale structure; features near the diagonal are short-lived and are often treated as less significant in applications.  The equivalent interval representation is the \emph{persistence barcode}, in which each class is drawn as an interval $[b_i,d_i)$.

In practice, we pass the graph-derived matrix $\mathbf{D}$ to Ripser as a precomputed condensed matrix and compute $\Dgm_0$ and $\Dgm_1$.  In the present time-series setting, these two diagrams have natural interpretations:
\begin{itemize}
  \item \textbf{$\Dgm_0$ (connected components).} At $\varepsilon=0$ every vertex is isolated.  As $\varepsilon$ grows, edges appear and components merge.  A short-lived $H_0$ class therefore indicates that the corresponding components are already close according to the chosen graph-based matrix.
  \item \textbf{$\Dgm_1$ (loops).} A class in $H_1$ is born when edges create a cycle that is not yet filled by triangles, and it dies when later simplices fill that loop.  For time-series networks, these features reflect cyclic or recurrent dynamics: a periodic signal typically yields one dominant loop, whereas a chaotic signal often yields multiple loops with a broader range of persistence values~\cite{Myers2019_PRE}.
\end{itemize}

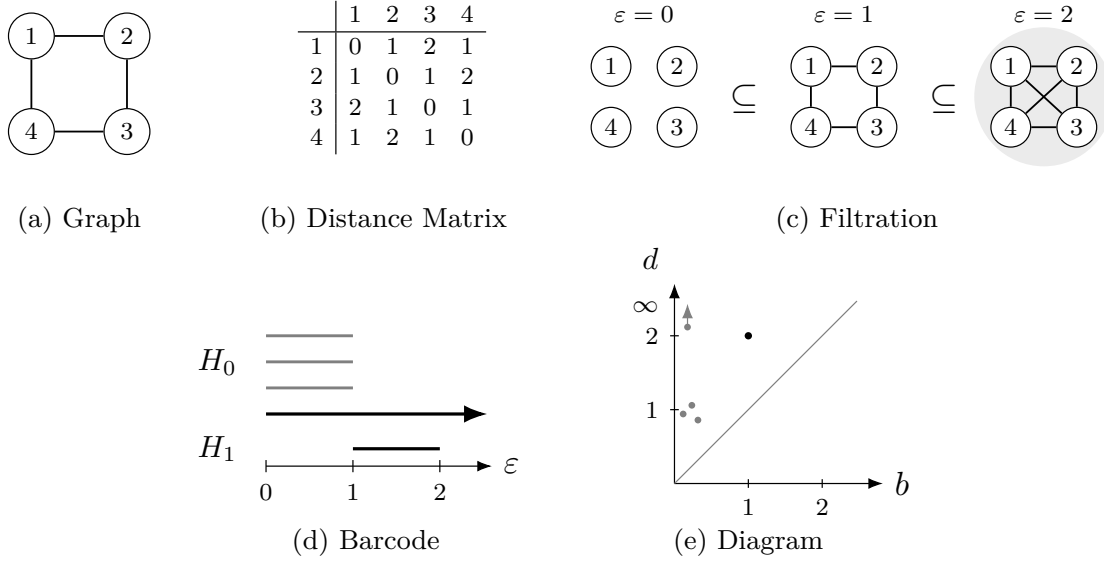
\begin{figure}[t]
\centering
\input{figures/rips-toy-filtration}
\caption{A single toy example summarizing the Section~1.2 route from graph to persistence.  Panel (a) shows an unweighted 4-cycle; panel (b) shows the shortest-path distance matrix derived from that graph; panel (c) shows the filtration at $\varepsilon=0,1,2$, where the complex has four connected components, then one loop, and then a filled loop; panels (d) and (e) show the corresponding barcode and persistence diagram.}
\label{fig:rips-toy-filtration}
\end{figure}

\Cref{fig:rips-toy-filtration} assembles this route into a single toy example: an unweighted 4-cycle, its shortest-path matrix, the resulting thresholded complexes, and the corresponding persistence barcode and diagram.

This topological summary is useful only if it changes in a controlled way under small perturbations of the filtration values.  The next subsection states the stability result that makes this precise.

\subsection{The Stability Theorem}
\label{sec:stability}

The remaining theoretical ingredient is robustness: once the filtration is built, small changes in the filtration values should produce only small changes in the persistence diagram.  In the present setting, a filtration function assigns to each simplex the threshold at which it enters the flag filtration, with the monotonicity condition that if $\sigma\subseteq\tau$, then $f(\sigma)\leq f(\tau)$.  Thus, one can think of two filtrations $f$ and $g$ as two nearby ways of assigning entry scales to the same underlying simplicial complex.

\begin{theorem}[Stability~{\cite{CohenSteiner2007}}]\label{thm:stability}
Let $f$ and $g$ be two filtration functions on the same finite simplicial complex~$K$.  Then, for every homological dimension~$p$,
\begin{equation}\label{eq:stability}
  \dB\bigl(\Dgm_p(f),\, \Dgm_p(g)\bigr) \;\leq\; \|f - g\|_\infty,
\end{equation}
where $\dB$ denotes the \emph{bottleneck distance}, i.e., the largest displacement needed in an optimal matching between the two diagrams.
\end{theorem}

In words, if every simplex enters at nearly the same scale in the two filtrations, then the resulting persistence diagrams cannot differ by more than that worst-case discrepancy.  In the present setting, the theorem applies after the graph-based matrix has been converted into simplex entry values: small perturbations in those values lead to small perturbations in $\Dgm_p$.

\subsection{Software and Workflow}
\label{sec:tools}

The software ecosystem for topological data analysis is broad, with mature libraries in several programming languages covering persistent-homology computation, Mapper-style visualisation, diagram vectorisation, and end-to-end machine-learning workflows.  A complete survey is beyond our scope.  In this study, the implementation is carried out primarily in Python and centres on the teaspoon library~\cite{Khasawneh2025}, which is designed for \emph{topological signal processing} and supports the parameter-selection procedures detailed later in \cref{sec:takens,sec:delay-selection,sec:dimension-selection}, the transition- and proximity-based network constructions of \cref{sec:networks}, persistence computation via a Ripser backend~\cite{Bauer2021}, and downstream benchmarking and machine-learning utilities.  For the visibility-based constructions, we use the dedicated \texttt{ts2vg} package to generate the NVG and HVG efficiently.  This combination is a better fit here than a persistence-only library because the present pipeline depends as much on upstream representation design as on persistence computation itself.

At the workflow level, the implementation proceeds as follows: (i)~construct the graph representation of the series, including delay-coordinate reconstruction and parameter selection when required (\cref{sec:takens,sec:delay-selection,sec:dimension-selection,sec:networks}); (ii)~compute the graph-derived matrix $\mathbf{D}$ (\cref{sec:distance-matrix}); (iii)~pass $\mathbf{D}$ to the Ripser backend for Vietoris--Rips persistence; and (iv)~vectorise the resulting diagrams (\cref{sec:vectorization}).

With the downstream route of this section held fixed, the remaining design question is upstream: how should the raw signal be converted into a graph in the first place?  We therefore turn next to the five network constructions that feed the common route.

\section{Network Construction Methods}
\label{sec:networks}

This section introduces the five network construction methods used in the chapter.  Each takes a univariate time series $x = [x_1, x_2, \ldots, x_L]$ as input and produces a graph representation.  Visibility and proximity methods yield undirected graphs directly, whereas transition methods such as the OPN and CGSSN naturally produce directed weighted networks that are symmetrized before entering the standardized route of \cref{sec:common-route}.  Thus the methodological variation lies here, at the graph-construction stage; once the graph has been built, the downstream pipeline is held fixed.

We begin with a brief overview of the three families and five methods.  The five methods fall into three families.  Visibility methods operate directly on the raw signal and connect time stamps by geometric line-of-sight criteria.  They are parameter-free and produce connected graphs, with their specific invariance properties discussed in \cref{sec:visibility}~\cite{Lacasa2008}.  Transition methods first construct a delay embedding and then convert each embedded point into a symbolic state, yielding compact graphs at the cost of parameter choices.  Proximity methods also begin with a delay embedding, but they preserve the reconstructed geometry more directly by connecting nearby states in phase space~\cite{Donner2010}.  \Cref{tab:comparison} previews the main tradeoffs.

\subsection{Takens' Delay Embedding}
\label{sec:takens}

Three of the five methods begin by lifting the scalar series into a higher-dimensional state space.  The ordinal partition network (OPN), the coarse-grained state-space network (CGSSN), and the $k$-nearest neighbor ($k$-NN) graph all rely on this reconstructed representation.

\begin{definition}[Delay embedding]\label{def:delay-embedding}
Given a uniformly sampled time series $x = [x_1, \ldots, x_L]$, an embedding dimension $n \geq 2$, and a time delay $\tau \geq 1$, the \emph{delay embedding} maps $x$ to a sequence of vectors in $\R^n$:
\begin{equation}\label{eq:delay-embedding}
  X_i = \bigl(x_i,\; x_{i+\tau},\; x_{i+2\tau},\; \ldots,\; x_{i+(n-1)\tau}\bigr), \qquad i = 1, \ldots, N,
\end{equation}
where $N = L - (n{-}1)\tau$.  The collection $\chi = \{X_1, \ldots, X_N\} \subset \R^n$ is the \emph{reconstructed point cloud}.
\end{definition}

Classically, Takens' theorem provides the theoretical justification for this construction: for generic smooth dynamics and observables, the delay-coordinate map is an embedding once $n \geq 2d + 1$, where $d$ is the manifold dimension~\cite{Takens1981}.  In practice, however, the attractor dimension is usually unknown, and finite noisy data rarely fit the theorem's ideal assumptions exactly.  For that reason, the embedding parameters $n$ and $\tau$ are typically chosen by data-driven heuristics rather than by direct application of the theorem.

\subsubsection{Delay Selection ($\tau$)}
\label{sec:delay-selection}

Several data-driven heuristics are used in practice to choose~$\tau$.  We emphasize two complementary approaches.

\paragraph{Mutual information (MI).}
A standard information-theoretic approach chooses a delay that reduces redundancy between $x_i$ and its delayed counterpart $x_{i+\tau}$ without making successive coordinates completely unrelated.  The mutual information function is
\begin{equation}\label{eq:mi}
  I(\tau) = \sum_{x_i,\, x_{i+\tau}} p\bigl(x_i,\, x_{i+\tau}\bigr) \ln \frac{p\bigl(x_i,\, x_{i+\tau}\bigr)}{p(x_i)\,p(x_{i+\tau})},
\end{equation}
where $p$ denotes the corresponding empirical marginal and joint distributions.  In practice, a common heuristic selects the \emph{first local minimum} of $I(\tau)$, which marks an early delay at which adjacent embedding coordinates are appreciably less redundant while still reflecting the underlying dynamics.

\paragraph{Multiscale permutation entropy (MsPE).}
For graph-based embeddings such as the OPN, one can instead examine a permutation-entropy curve across candidate delays~\cite{AzizArif2005,Myers2024_delay}.  For a candidate delay~$\tau$, MsPE constructs delay vectors and maps each to its ordinal pattern $\pi_i \in S_n$.  The relative frequencies define the \emph{permutation entropy at scale~$\tau$}:
\begin{equation}\label{eq:mspe}
  H(\tau) = -\sum_{\pi \in S_n} p(\pi;\tau) \log_2 p(\pi;\tau),
\end{equation}
with normalization $h(\tau) = H(\tau)/\log_2(n!) \in [0,1]$.  At small delays, successive embedding coordinates are nearly identical, so the trajectory visits only a few ordinal patterns and the entropy is low.  As $\tau$ increases, entropy typically rises.  For periodic or quasiperiodic signals, one often observes a \emph{resonance dip} associated with the dominant period, and that feature can serve as a practical cue for choosing~$\tau$~\cite{Myers2024_delay}.  For chaotic signals where no clear resonance dip appears, one may instead use the first prominent entropy maximum as a fallback, aiming for a delay that yields a richer spread of ordinal patterns.

\subsubsection{Dimension Selection ($n$)}
\label{sec:dimension-selection}

In practice, the embedding dimension is usually chosen heuristically, most commonly via the false-nearest-neighbors test.

\paragraph{False nearest neighbors (FNN).}
The false-nearest-neighbors test of Kennel et~al.~\cite{Kennel1992} iteratively embeds the time series in progressively higher dimensions.  If the current dimension~$n$ is too low, the projection forces distant points in the true phase space to falsely appear as close neighbors.  For each embedded point $X_i^{(n)}$ with nearest neighbor $X_j^{(n)}$ at distance $R_i^{(n)} = \|X_i^{(n)} - X_j^{(n)}\|$, a neighbor is declared \emph{false} if
\begin{equation}\label{eq:fnn}
  \frac{|x_{i+n\tau} - x_{j+n\tau}|}{R_i^{(n)}} > R_{\mathrm{tol}},
\end{equation}
where $R_{\mathrm{tol}}$ is a threshold (typically $R_{\mathrm{tol}} = 15$).  The embedding dimension is the smallest~$n$ for which the FNN fraction drops below $1$--$2\%$.  Small values such as $n = 3$ to $7$ are common in practice.

For ordinal methods such as the OPN, one also keeps $n$ modest in practice because the number of possible states grows rapidly with~$n$, which can make the transition counts sparse.  In this setting, the MsPE framework naturally extends from delay to dimension selection: one sweeps $n$ over the candidate range $[n_{\min}, n_{\max}]$ and selects the value that maximises the normalised permutation entropy, ensuring a rich spread of visited ordinal patterns.  This approach provides a permutation-native alternative to FNN for the OPN.

\begin{remark}[Practical starting values]\label{rem:parameter-ranges}
In the examples below, we treat $\tau \in [1,50]$ and $n \in [3,7]$ as practical search ranges for the embedding parameters, informed by the heuristic discussions above and by recent delay-selection work of Myers et~al.~\cite{Myers2024_delay}.  For the CGSSN, the bin count $b$ is application-dependent and is tuned to balance state-space resolution against sparsity~\cite{WangTian2016,Myers2023_CGSSN}.  For the $k$-NN graph, $k$ is adjusted to maintain connectivity without overly densifying the graph.  For diffusion distance (see \cref{sec:distance-matrix}), Myers et~al.~\cite{Myers2023_SIAM} used the heuristic $d < t < 3d$, where $d$ denotes the graph diameter.  These are tuning guidelines rather than universal defaults.
\end{remark}

\begin{remark}[Non-uniform sampling]\label{rem:nonuniform}
\Cref{def:delay-embedding} assumes uniformly sampled data.  For irregularly sampled time series, a common practical step is to resample or interpolate the observations onto a uniform grid before applying the standard delay embedding used here.  The interpolation method and target sampling rate are application-dependent and can materially affect the reconstructed attractor geometry.
\end{remark}

\subsection{Visibility Graphs}
\label{sec:visibility}

Visibility graphs~\cite{Lacasa2008} operate directly on the sampled signal: each observation $(t_i, y_i)$ is treated as a vertical bar, and two bars are connected when they can ``see'' each other.

\subsubsection{Natural Visibility Graph (NVG)}

\begin{definition}[Natural visibility criterion]\label{def:nvg}
Nodes $v_i$ and $v_j$ ($t_i < t_j$) are connected if and only if every intermediate observation $(t_k, y_k)$ with $t_i < t_k < t_j$ satisfies
\begin{equation}\label{eq:nvg}
  y_k < y_j + (y_i - y_j)\,\frac{t_j - t_k}{t_j - t_i}.
\end{equation}
\end{definition}

The resulting graph is always connected and undirected, and it is invariant under affine transformations of the signal values.  Local maxima often become hubs, and the degree distribution can reflect properties of the generating process: periodic series tend to produce regular graphs, random series often yield exponential behavior, and fractal series may exhibit scale-free structure.

\subsubsection{Horizontal Visibility Graph (HVG)}

\begin{definition}[Horizontal visibility criterion]\label{def:hvg}
Nodes $v_i$ and $v_j$ ($t_i < t_j$) are connected if and only if $y_k < \min(y_i, y_j)$ for all $t_i < t_k < t_j$.
\end{definition}

Every HVG edge is also an NVG edge, so the HVG is a subgraph of the NVG.  Luque et~al.~\cite{Luque2009} derived exact results for the degree distribution under i.i.d.\ random variables, providing a null model for hypothesis testing.

Both variants are parameter-free, which is a significant practical advantage.  The na\"{\i}ve NVG algorithm has $O(L^2)$ complexity, though divide-and-conquer implementations achieve $O(L \log L)$.  In either case, the graph has $L$ nodes, one per time step.

\begin{remark}\label{rem:visibility-limitation}
Series with different amplitude ranges but identical ordering structure can produce the same unweighted graph.  This insensitivity to amplitude can be either an advantage (robustness to scaling) or a disadvantage (loss of discriminative information), depending on the task.
\end{remark}

\Cref{ex:nvg} applies both criteria to the same five-point signal.  The signal is designed so that the global maximum at $t=4$ becomes an obvious hub---illustrating the claim that local maxima attract long-range edges---while one node pair falls exactly in the gap between the NVG and HVG rules, making the difference between the two variants concrete at the smallest possible scale.

\begin{example}\label{ex:nvg}
Consider $x = [1, 2, 1, 5, 2]$ with $t = [1,2,3,4,5]$.

\textbf{Step~1.}  Each observation $(t_i, x_i)$ becomes a node.  All consecutive pairs $(i, i{+}1)$ are automatically connected in both variants.

\textbf{Step~2 (NVG).}  For non-adjacent pairs, check whether the joining line clears every intermediate bar.  Two long-range edges pass:
\begin{itemize}
  \item $(1,4)$: the line from $(1,1)$ to $(4,5)$ gives heights $2.\overline{3}$ and $3.\overline{6}$ at $t=2$ and $t=3$; since $x_2=2<2.\overline{3}$ and $x_3=1<3.\overline{6}$, the edge exists.
  \item $(2,4)$: line height at $t=3$ is $3.5$; $x_3=1<3.5$.  Edge exists.
\end{itemize}
Node~$4$ ($x_4=5$, global maximum) is the hub, with degree~$4$.

\textbf{Step~3 (HVG).}  Replace the line test with the stricter threshold $x_k < \min(x_i, x_j)$:
\begin{itemize}
  \item $(1,4)$: $\min(x_1,x_4)=\min(1,5)=1$, but $x_2=2\geq 1$.  \emph{Edge blocked.}
  \item $(2,4)$: $\min(2,5)=2$; $x_3=1<2$.  Edge retained.
\end{itemize}
The NVG has $6$ edges and the HVG has $5$; the sole difference is edge $(1,4)$.  \Cref{fig:visibility-toy} illustrates this comparison.
\end{example}

\begin{center}
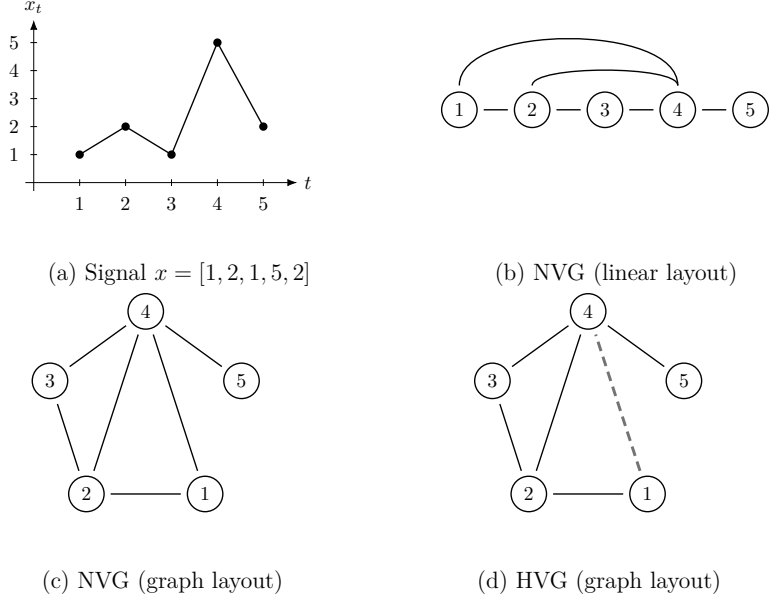

\input{figures/section13-visibility-toy}
\captionof{figure}{Visibility-graph comparison for $x=[1,2,1,5,2]$.  Panel~(a): signal; node~$4$ ($x_4=5$) is the global maximum.  Panel~(b): NVG in linear layout; long-range edges $(1,4)$ and $(2,4)$ appear as nested arcs.  Panel~(c): NVG redrawn in a graph layout with no spatial constraints; node~$4$ is the hub of degree~$4$.  Panel~(d): HVG in the same graph layout; the dashed grey edge marks the absent edge $(1,4)$, blocked because $x_2=2\geq\min(x_1,x_4)=1$.}
\label{fig:visibility-toy}
\end{center}

\subsection{Ordinal Partition Networks}
\label{sec:opn}

The ordinal partition network (OPN) converts a time series into a transition network whose nodes represent ordinal patterns.  This construction was introduced by McCullough et~al.~\cite{McCullough2015} and is rooted in the permutation entropy framework of Bandt and Pompe~\cite{BandtPompe2002}.

\begin{definition}[Ordinal partition network]\label{def:opn}
Given the delay embedding $\chi = \{X_1, \ldots, X_N\}$ with dimension~$n$ and delay~$\tau$, each vector $X_i$ is mapped to its \emph{ordinal pattern}: the permutation $\pi_i \in S_n$ that sorts the coordinates into ascending order.  When ties occur, they are broken by temporal order, so earlier coordinates precede later ones.  The OPN is the directed, weighted graph $G_{\mathrm{OPN}} = (V, E, \omega)$ where:
\begin{itemize}
  \item $V = \{\pi \in S_n \mid \pi = \pi_i \text{ for some } i\}$ (observed permutations),
  \item $(u, v) \in E$ if the series transitions from pattern~$u$ to pattern~$v$ at some consecutive pair $(i, i{+}1)$,
  \item $\omega(u, v)$ counts the number of such transitions.
\end{itemize}
The graph is symmetrized as $\mathbf{A}_{\mathrm{sym}} = \mathbf{A} + \mathbf{A}^T$, yielding an undirected weighted graph where $\omega_{\mathrm{sym}}(u,v) = \omega(u,v) + \omega(v,u)$ counts total transitions in both directions.
\end{definition}

Because the ordinal pattern records only the \emph{relative ordering} of values, the OPN is inherently robust to monotone amplitude transformations.  This invariance is useful when amplitude scaling is not itself the object of interest.  In canonical examples, periodic dynamics produce cycle-like networks, whereas chaotic dynamics generate more richly connected structures~\cite{McCullough2015}.  The tradeoff is that amplitude information is discarded entirely: two signals with identical ordinal sequences but different amplitude ranges yield the same OPN.

To see concretely how ordinal patterns are assigned and how repeated patterns accumulate transition weights, consider \cref{ex:opn}.  The signal is long enough to revisit two of the five observed patterns---creating multi-weight edges after symmetrization---without being so long that the bookkeeping becomes unwieldy.

\begin{example}\label{ex:opn}
Consider $x = [1,3,2,5,4,2,4,1,3]$ with $n=3$, $\tau=1$, giving $N = 7$ delay vectors.

\textbf{Step~1.}  Form the delay vectors:
\[
  X_1=[1,3,2],\; X_2=[3,2,5],\; X_3=[2,5,4],\; X_4=[5,4,2],\; X_5=[4,2,4],\; X_6=[2,4,1],\; X_7=[4,1,3].
\]

\textbf{Step~2.}  Map each vector to its ordinal pattern~$\pi$ (indices sorted by ascending value; ties broken by earlier index):
\begin{align*}
  X_1=[1,3,2] &\;\to\; \pi_1=(0,2,1), \quad\text{since } 1<2<3; \\
  X_2=[3,2,5] &\;\to\; \pi_2=(1,0,2), \quad\text{since } 2<3<5; \\
  X_3=[2,5,4] &\;\to\; \pi_1=(0,2,1), \quad\text{same pattern as }X_1; \\
  X_4=[5,4,2] &\;\to\; \pi_5=(2,1,0), \quad\text{since } 2<4<5; \\
  X_5=[4,2,4] &\;\to\; \pi_2=(1,0,2), \quad\text{since } 2<4{=}4\text{ (tie broken: index~$0<$index~$2$)}; \\
  X_6=[2,4,1] &\;\to\; \pi_4=(2,0,1), \quad\text{since } 1<2<4; \\
  X_7=[4,1,3] &\;\to\; \pi_3=(1,2,0), \quad\text{since } 1<3<4.
\end{align*}
Five distinct patterns are observed: $\pi_1,\pi_2,\pi_3,\pi_4,\pi_5$ become the five nodes.

\textbf{Step~3.}  Each consecutive pair $(X_i, X_{i+1})$ gives a directed transition:
\[
  \pi_1\!\to\!\pi_2,\quad
  \pi_2\!\to\!\pi_1,\quad
  \pi_1\!\to\!\pi_5,\quad
  \pi_5\!\to\!\pi_2,\quad
  \pi_2\!\to\!\pi_4,\quad
  \pi_4\!\to\!\pi_3.
\]

\textbf{Step~4 (symmetrize).}  Merge directed counts into undirected weights $\omega_{\mathrm{sym}}(u,v)=\omega(u,v)+\omega(v,u)$:
\begin{itemize}
  \item $\{\pi_1,\pi_2\}$: each direction appears once $\Rightarrow$ weight~$2$.
  \item $\{\pi_1,\pi_5\}$, $\{\pi_2,\pi_4\}$, $\{\pi_2,\pi_5\}$, $\{\pi_3,\pi_4\}$: one direction each $\Rightarrow$ weight~$1$.
\end{itemize}
Node $\pi_2$ is the hub of degree~$3$; see \cref{fig:opn-toy}.
\end{example}

\noindent\makebox[\textwidth][c]{%
  \input{figures/section13-opn-toy}%
}
\captionof{figure}{Ordinal partition network for $x=[1,3,2,5,4,2,4,1,3]$ with $n=3$, $\tau=1$.
  Panel~(a): signal; the shaded region marks the first sliding window of width~$n=3$.
  Panel~(b): OPN with five observed pattern-nodes $\pi_1,\ldots,\pi_5$
  (where $\pi_1{=}(0,2,1)$, $\pi_2{=}(1,0,2)$, $\pi_3{=}(1,2,0)$, $\pi_4{=}(2,0,1)$, $\pi_5{=}(2,1,0)$).
  The bold edge $\{\pi_1,\pi_2\}$ carries weight~$2$ because both directions
  $\pi_1\!\to\!\pi_2$ and $\pi_2\!\to\!\pi_1$ each appear once;
  all other edges carry weight~$1$.  Node $\pi_2$ is the hub of degree~$3$.
}
\label{fig:opn-toy}

\subsection{Coarse-Grained State-Space Networks}
\label{sec:cgssn}

The coarse-grained state-space network (CGSSN) used here follows the phase-space coarse-graining construction of Wang and Tian~\cite{WangTian2016}.  Like the OPN, it begins from a delay embedding and records temporal transitions, but it replaces ordinal patterns with discretized amplitude states.  It is also related to earlier phase-space reconstruction-network approaches such as Gao and Jin~\cite{GaoJin2009}.

\begin{definition}[Coarse-grained state-space network]\label{def:cgssn}
Given the delay embedding $\chi$ with dimension~$n$ and delay~$\tau$, and a bin count $b \geq 2$, the amplitude range $[\min(x), \max(x)]$ is partitioned into $b$ equal-width bins.  Each coordinate of $X_i$ is mapped to its bin index $\rho_i(j) \in \{0, \ldots, b{-}1\}$, and the state is encoded as
\begin{equation}\label{eq:cgssn-state}
  s_i = 1 + \sum_{j=0}^{n-1} \rho_i(j) \cdot b^j.
\end{equation}
The CGSSN is constructed analogously to the OPN: nodes are observed states, directed edges record temporal transitions, and weights count transition frequencies.  The graph is symmetrized for the persistence pipeline.
\end{definition}

The total number of possible states is $b^n$, which can be much larger than the $n!$ states of the OPN; for instance, $b=8$ and $n=4$ gives $4096$ states versus $24$ for the OPN.  This larger state space increases descriptive power but also raises the risk of sparsity when the time series is short.  Myers et~al.~\cite{Myers2023_CGSSN} studied persistent homology on CGSSNs for dynamic-state detection, including comparisons with OPN-based constructions.

\Cref{ex:cgssn} below uses the coarsest possible binning ($b=2$) to keep the state count small while demonstrating the key steps: how the amplitude range is partitioned into bins, how the state-index formula encodes a multi-dimensional bin tuple into a single integer, and how two delay vectors that land in the same bin tuple produce a repeated state that yields a weighted edge after symmetrization.

\begin{example}\label{ex:cgssn}
Consider $x=[1,4,7,2,5,8,3,6,1]$ with $n=3$, $\tau=1$, $b=2$, giving $N=7$ delay vectors.

\textbf{Step~1.}  Form the delay vectors:
\[
  X_1=[1,4,7],\; X_2=[4,7,2],\; X_3=[7,2,5],\; X_4=[2,5,8],\; X_5=[5,8,3],\; X_6=[8,3,6],\; X_7=[3,6,1].
\]

\textbf{Step~2.}  Partition $[\min x, \max x]=[1,8]$ into $b=2$ equal bins of width $3.5$:
\[
  \text{bin~}0 = [1,\, 4.5), \qquad \text{bin~}1 = [4.5,\, 8].
\]
Assign each coordinate to its bin and compute $s_i = 1 + \sum_{j=0}^{2}\rho_i(j)\cdot 2^j$:
\begin{align*}
  X_1=[1,4,7]: &\; \rho=[0,0,1],\quad s_1 = 1+0+0+4 = 5. \\
  X_2=[4,7,2]: &\; \rho=[0,1,0],\quad s_2 = 1+0+2+0 = 3. \\
  X_3=[7,2,5]: &\; \rho=[1,0,1],\quad s_3 = 1+1+0+4 = 6. \\
  X_4=[2,5,8]: &\; \rho=[0,1,1],\quad s_4 = 1+0+2+4 = 7. \\
  X_5=[5,8,3]: &\; \rho=[1,1,0],\quad s_5 = 1+1+2+0 = 4. \\
  X_6=[8,3,6]: &\; \rho=[1,0,1],\quad s_6 = 1+1+0+4 = 6. \\
  X_7=[3,6,1]: &\; \rho=[0,1,0],\quad s_7 = 1+0+2+0 = 3.
\end{align*}
Vectors $X_3$ and $X_6$ share state~$6$; $X_2$ and $X_7$ share state~$3$.  Five distinct states $\{3,4,5,6,7\}$ become the five nodes.

\textbf{Step~3.}  State sequence $5\to3\to6\to7\to4\to6\to3$ gives directed transitions:
\[
  5\!\to\!3,\quad 3\!\to\!6,\quad 6\!\to\!7,\quad 7\!\to\!4,\quad 4\!\to\!6,\quad 6\!\to\!3.
\]

\textbf{Step~4 (symmetrize).}  Merge directed counts into undirected weights:
\begin{itemize}
  \item $\{3,6\}$: directions $3\!\to\!6$ and $6\!\to\!3$ each appear once $\Rightarrow$ weight~$2$.
  \item $\{3,5\}$, $\{6,7\}$, $\{4,7\}$, $\{4,6\}$: one direction each $\Rightarrow$ weight~$1$.
\end{itemize}
State~$6$ is the hub of degree~$3$; see \cref{fig:cgssn-toy}.
\end{example}

\noindent\makebox[\textwidth][c]{%
  \input{figures/section13-cgssn-toy}%
}

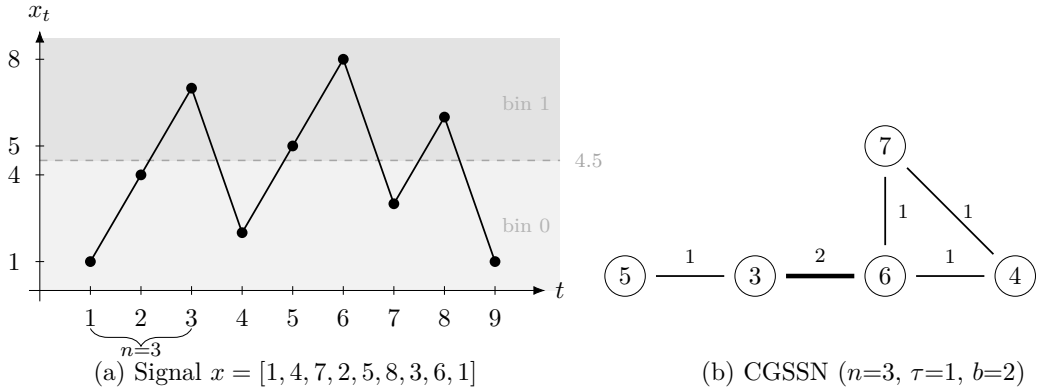
\captionof{figure}{Coarse-grained state-space network for $x=[1,4,7,2,5,8,3,6,1]$ with
  $n=3$, $\tau=1$, $b=2$.
  Panel~(a): signal; the shaded bands mark bin~$0=[1,4.5)$ (light) and bin~$1=[4.5,8]$ (medium),
  with the dashed boundary line at~$4.5$.
  Panel~(b): CGSSN with five state-nodes.
  The bold edge $\{3,6\}$ carries weight~$2$ because both directions appear once;
  all other edges carry weight~$1$.  State~$6$ is the hub of degree~$3$.
}
\label{fig:cgssn-toy}

\subsection{\texorpdfstring{$k$}{k}-Nearest Neighbor Graphs}
\label{sec:knn}

Among the methods considered here, the $k$-nearest neighbor graph preserves the reconstructed geometry most directly.  Rather than symbolizing the embedded vectors, it connects points that are nearby in the delay-coordinate space.

\begin{definition}[$k$-nearest neighbor graph]\label{def:knn}
Given the delay embedding $\chi = \{X_1, \ldots, X_N\} \subset \R^n$ and a neighborhood parameter $k \geq 1$, the $k$-NN graph $G_{k\text{-NN}} = (V, E)$ has vertex set
\[
  V = \{v_1, \ldots, v_N\}, \qquad v_i \leftrightarrow X_i,
\]
and edge set
\[
  E = \bigl\{(v_i, v_j) \mid X_j \in \NN_k(X_i) \text{ or } X_i \in \NN_k(X_j)\bigr\},
\]
where $\NN_k(X_i)$ denotes the $k$ nearest neighbors of $X_i$ under the Euclidean metric.
\end{definition}

The ``or'' condition in the edge set produces an \emph{undirected} graph, namely the symmetric union of directed $k$-NN neighborhoods, which is required for the subsequent symmetric distance matrix.  The number of edges is at most $kN$, though typically fewer due to reciprocal neighbors.  The neighborhood size~$k$ controls graph density: too small and the graph may disconnect; too large and local structure is lost.  A common starting point is to take $k$ on the order of $\sqrt{N}$ and then adjust as needed to maintain connectivity without overly densifying the graph.

The $k$-NN construction is closely related to manifold learning: in the Isomap algorithm~\cite{Tenenbaum2000}, shortest-path distances on a $k$-NN graph approximate geodesic distances on the underlying manifold.  When we compute the distance matrix on the $k$-NN graph (\cref{sec:distance-matrix}), we obtain a discrete approximation to the intrinsic geometry of the reconstructed attractor.

\begin{remark}\label{rem:knn-vs-transition}
The $k$-NN graph preserves metric information that the OPN discards (amplitude) and that the CGSSN only coarsely approximates.  However, it has $N = L - (n{-}1)\tau$ nodes---potentially much larger than the $n!$ or $b^n$ nodes of the transition networks---making the Rips computation more expensive.
\end{remark}

\Cref{ex:knn} is designed to highlight precisely this distinction.  The signal contains a repeated amplitude value so that two delay vectors lie close together in $\R^3$, and the hub node that emerges is determined entirely by metric proximity in delay-coordinate space rather than by temporal succession---a structural feature that neither the OPN nor the CGSSN would reveal.

\begin{example}\label{ex:knn}
Consider $x=[1,2,1,5,2,3,1]$ with $n=3$, $\tau=1$, $k=2$, giving $N=5$ delay vectors.

\textbf{Step~1.}  Form the delay vectors:
\[
  X_1=[1,2,1],\quad X_2=[2,1,5],\quad X_3=[1,5,2],\quad
  X_4=[5,2,3],\quad X_5=[2,3,1].
\]

\textbf{Step~2.}  Compute all pairwise Euclidean distances in $\R^3$, listed in increasing order:
\[
  d(X_1,X_5)=\!\sqrt{2},\quad d(X_3,X_5)=\!\sqrt{6},\quad d(X_1,X_3)=\!\sqrt{10},\quad
  d(X_2,X_4)=d(X_4,X_5)=\!\sqrt{14},
\]
\[
  d(X_1,X_2)=\!\sqrt{18},\quad d(X_1,X_4)=d(X_2,X_5)=\!\sqrt{20},\quad
  d(X_2,X_3)=d(X_3,X_4)=\!\sqrt{26}.
\]

\textbf{Step~3.}  Find the $2$ nearest neighbors of each point:
\begin{align*}
  &\NN_2(X_1)=\{X_5,X_3\},\quad
   \NN_2(X_2)=\{X_4,X_1\},\quad
   \NN_2(X_3)=\{X_5,X_1\},\\
  &\NN_2(X_4)=\{X_2,X_5\},\quad
   \NN_2(X_5)=\{X_1,X_3\}.
\end{align*}

\textbf{Step~4 (symmetric union).}  Connect $v_i$ and $v_j$ whenever $X_j\in\NN_k(X_i)$ or $X_i\in\NN_k(X_j)$:
\[
  E = \bigl\{(1,2),(1,3),(1,5),(2,4),(3,5),(4,5)\bigr\}.
\]
The graph has $5$ nodes and $6$ edges.  Node~$v_1$ (corresponding to $X_1=[1,2,1]$) is a hub of degree~$3$; see \cref{fig:knn-toy}.
\end{example}

\noindent\makebox[\textwidth][c]{\input{figures/section13-knn-toy}}

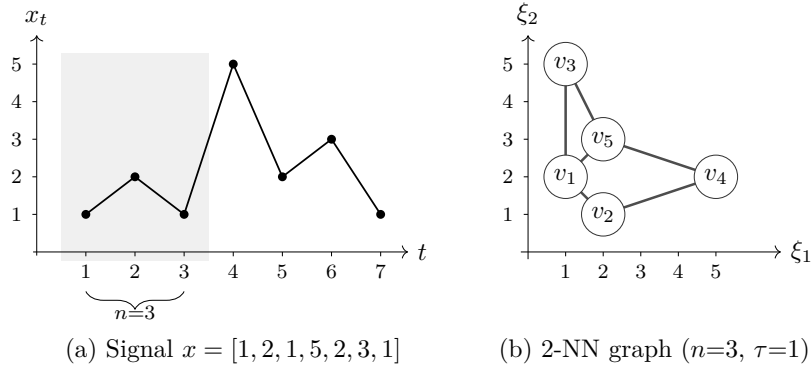
\captionof{figure}{$k$-nearest neighbor graph for $x=[1,2,1,5,2,3,1]$ with $n=3$, $\tau=1$, $k=2$.  Panel~(a): signal; the shaded region marks the first sliding window of width $n=3$.  Panel~(b): $2$-NN graph; nodes placed at the first two delay-vector components $(\xi_1,\xi_2)$ to reveal the geometric structure that drives edge formation.  Edges connect pairs whose 3-D Euclidean distance is small enough to fall within each other's 2-nearest-neighbor list.  Node $v_1$ is the hub of degree~$3$.}
\label{fig:knn-toy}

\subsection{Relationships Between the Methods}
\label{sec:relationships}

Following the taxonomy surveyed by Silva et~al.~\cite{Silva2021}, the five methods represent three distinct ways of turning a time series into a network.  Visibility graphs (NVG and HVG) work directly on the observed signal and emphasize geometric structure in time; HVG imposes a stricter horizontal threshold and is a subgraph of NVG.  The OPN and CGSSN first reconstruct a state space and then record temporal movement through that space as transitions between symbolic states.  The $k$-NN graph also begins with a delay embedding, but instead of encoding temporal succession it connects states that are nearby in the reconstructed geometry.  Thus the transition methods emphasize symbolic evolution, whereas the $k$-NN graph emphasizes spatial proximity.

These differences matter for both graph size and the kind of information preserved.  Visibility graphs have one node per observation and therefore remain close to the raw series.  The OPN is more compact, with at most $n!$ nodes, because it retains only ordinal information; this makes it insensitive to monotone rescaling but discards amplitude.  The CGSSN retains coarse amplitude information by replacing ordinal states with discretized bins, at the cost of a larger state space of size at most $b^n$.  The $k$-NN graph preserves metric structure most directly, but it inherits the full embedded sample size $N = L - (n{-}1)\tau$ and is therefore often the most expensive option downstream.

The constructions also suggest different sensitivities to perturbation.  Because the OPN depends only on relative ordering, it is less affected by monotone amplitude transformations.  The CGSSN can react more strongly to small amplitude changes when those changes shift points across bin boundaries.  Visibility graphs may gain or lose edges when local peaks alter line-of-sight relations.  The $k$-NN graph reflects noise in the reconstructed coordinates directly, although increasing~$k$ can smooth small perturbations at the cost of weakening local detail.

\subsection{Comparative Summary}
\label{sec:method-comparison}

\begin{table}[ht]
\centering\small
\caption{Comparative summary of the five network construction methods.}
\label{tab:comparison}
\begin{tabular}{@{}lllcl>{\raggedright\arraybackslash}p{3.4cm}@{}}
\toprule
\textbf{Method} & \textbf{Family} & \textbf{Encodes} & \textbf{Nodes} & \textbf{Parameters} & \textbf{Typical use case} \\
\midrule
NVG     & Visibility & Geometric structure    & $L$              & None                & More edges ($O(L\log L)$); baseline \\
HVG     & Visibility & Geometric structure    & $L$              & None                & Subgraph of NVG; stricter baseline \\
OPN     & Transition & Symbolic dynamics      & $\leq n!$        & $n$, $\tau$         & Compact symbolic representation \\
CGSSN   & Transition & Coarse state transitions & $\leq b^n$       & $n$, $\tau$, $b$    & Transition model retaining amplitude \\
$k$-NN  & Proximity  & Attractor geometry     & $L{-}(n{-}1)\tau$ & $n$, $\tau$, $k$    & Geometry-oriented representation \\
\bottomrule
\end{tabular}
\end{table}

This comparison closes the graph-construction part of the chapter.  With the upstream representation fixed and the downstream filtration already standardized, the only remaining methodological question is how to convert $\Dgm_0$ and $\Dgm_1$ into classifier-ready numerical features.

\section{From Persistence Diagrams to Feature Vectors}
\label{sec:vectorization}

Given the pair of persistence diagrams $\Dgm_0$ and $\Dgm_1$ produced by the common route of \cref{sec:common-route}, the next step is vectorization.  These diagrams are variable-size multisets and cannot be fed directly into standard classifiers, which expect fixed-length numerical vectors.  The goal of this section is to describe stable, finite-dimensional mappings from diagrams to $\R^d$.  Throughout, we reserve bold $\mathbf{D}$ for the graph-derived pairwise matrix of \cref{sec:distance-matrix} and use plain $D$ for persistence diagrams, to avoid notational conflict.

\subsection{Persistence Landscapes}
\label{sec:landscapes}

Persistence landscapes, introduced by Bubenik~\cite{Bubenik2015}, embed persistence diagrams into a Banach space of piecewise-linear functions.

Given a persistence diagram $D = \{(b_i, d_i)\}_{i=1}^{m}$, each birth--death pair defines a tent function
\begin{equation}\label{eq:tent}
  \Lambda_i(t) = \begin{cases}
    t - b_i     & \text{if } t \in \bigl(b_i,\, \tfrac{b_i + d_i}{2}\bigr], \\[3pt]
    d_i - t     & \text{if } t \in \bigl(\tfrac{b_i + d_i}{2},\, d_i\bigr), \\[3pt]
    0           & \text{otherwise.}
  \end{cases}
\end{equation}
Each tent peaks at the midpoint $(b_i + d_i)/2$ with height $(d_i - b_i)/2$.  The \emph{$k$-th persistence landscape} is
\begin{equation}\label{eq:landscape}
  \lambda_k(t) = \kmax\bigl\{\Lambda_i(t) \mid i = 1, \ldots, m\bigr\},
\end{equation}
where $\kmax$ denotes the $k$-th largest value.

\begin{theorem}[Landscape stability~{\cite{Bubenik2015}}]\label{thm:landscape-stability}
For any two persistence diagrams $D$ and $D'$,
\begin{equation}\label{eq:landscape-stability}
  |\lambda_k(t) - \lambda'_k(t)| \leq \dB(D, D')
\end{equation}
for all $k$ and $t$.
\end{theorem}

Persistence landscapes live in the Banach space $L^p(\mathbb{N} \times \R)$; for $p=2$ this becomes a Hilbert space, enabling means, variances, and inner products.  Bubenik~\cite{Bubenik2015} shows that the map from diagrams to their complete landscape functions is injective: two diagrams are equal if and only if all their landscape layers agree pointwise.  This guarantee applies to the full infinite-dimensional object; the finite discretization introduced in the next subsubsection is a deliberate lossy compression.

Composing \cref{thm:landscape-stability} with \cref{thm:stability} yields an end-to-end guarantee: $|\lambda_k(t) - \lambda'_k(t)| \leq \dB(D,D') \leq \|f - g\|_\infty$.  Small perturbations in the input filtration therefore produce only small changes in every landscape layer---the dominant topological signal is not washed out by noise.

\subsubsection{Discretization for Machine Learning}

For practical machine learning, the landscape functions must be compressed into finite vectors.  In the standardized pipeline adopted here, each landscape $\lambda_k$ is evaluated on a uniform grid of $M = 200$ points spanning $[\min(b_i), \max(d_i)]$.  We retain the first $J = 3$ landscape layers and compute their pointwise mean:
\begin{equation}\label{eq:mean-landscape}
  \bar{\lambda}(t_j) = \frac{1}{J}\sum_{k=1}^{J} \lambda_k(t_j),
\end{equation}
yielding a $200$-dimensional landscape feature vector per homology dimension.  These values are pragmatic chapter-level defaults rather than universal optima: a finer grid provides better resolution of the landscape shape while keeping the feature dimension tractable.  In a fully tuned predictive workflow, both $M$ and $J$ should be selected by validation or domain-specific constraints.

\subsection{Topological Summary Statistics}
\label{sec:summary-stats}

Scalar summary statistics complement the landscape representation by compressing global diagram properties into interpretable numbers.  Each statistic maps the entire diagram $D = \{(b_i, d_i)\}_{i \in I}$ to a single value; applied separately to $\Dgm_0$ and $\Dgm_1$, each statistic contributes two features to the final vector.

\paragraph{Persistent entropy.}
\begin{equation}\label{eq:entropy}
  E(D) = -\sum_{i \in I} \frac{\ell_i}{L_D} \log\!\Bigl(\frac{\ell_i}{L_D}\Bigr),
\end{equation}
where $\ell_i = d_i - b_i$ is the lifetime and $L_D = \sum_{i \in I} \ell_i$ is the total lifetime~\cite{Guzel2023}.  This measures the complexity of the lifetime distribution and is scale-invariant.

\paragraph{Amplitude (maximum lifetime).}
$A(D) = \max_{i \in I}(d_i - b_i)$, capturing the most prominent topological feature.

\paragraph{Total persistence.}
The total lifetime
\[
  L_D = \sum_{i \in I} (d_i - b_i),
\]
which measures the aggregate persistence mass of the diagram.  Under the standard diagram-metric convention, it is best interpreted as a lifetime sum rather than identified directly with a Wasserstein-$1$ distance to the empty diagram.

\paragraph{Number of points.}
The cardinality $|I|$, indicating the structural complexity of the underlying topology.

\paragraph{Adcock--Carlsson coordinates.}
Polynomial symmetric functions capturing interactions between birth times, death times, and lifetimes~\cite{AdcockCarlsson2016}:
\begin{equation}\label{eq:ac1}
  f_1(D) = \sum_{i \in I} b_i \,\ell_i, \qquad
  f_2(D) = \sum_{i \in I} (d_{\max} - d_i)\,\ell_i,
\end{equation}
\begin{equation}\label{eq:ac2}
  f_3(D) = \sum_{i \in I} b_i^2 \,\ell_i^4, \qquad
  f_4(D) = \sum_{i \in I} (d_{\max} - d_i)^2\,\ell_i^4,
\end{equation}
where $d_{\max} = \max_{i \in I} d_i$.

\paragraph{Landscape norm.}
The $L^1$ norm of the first landscape, $\|\lambda_1\|_1$, summarizes the dominant layer.  When tent functions do not overlap, a closed form is available~\cite{Akingbade2024}:
\begin{equation}\label{eq:landscape-norm}
  \|\lambda_1\|_1 = \frac{1}{4}\sum_i (d_i - b_i)^2.
\end{equation}

\subsection{Assembling the Final Feature Vector}
\label{sec:feature-assembly}

For each time series, the vectorization stage produces a feature vector by concatenating two blocks:
\begin{enumerate}
  \item \textbf{Landscape features.}  The mean persistence landscape for $H_0$ and $H_1$, each at $M = 200$ grid points, yielding $400$ features.
  \item \textbf{Summary statistics.}  The scalar summaries of \cref{sec:summary-stats} for both homology dimensions: persistent entropy, amplitude, total persistence, cardinality, four Adcock--Carlsson coordinates, and landscape norm, yielding $18$ features.
\end{enumerate}

The concatenated $418$-dimensional vector is passed to a standard classifier.  Fixing all vectorization choices before classification ensures that the feature vector has the same length for every time series in the dataset, so that classification accuracy reflects upstream design choices rather than variable feature geometry.

\subsection{Illustrative Example}
\label{sec:worked-example}

\Cref{fig:landscape-toy} illustrates the construction for the diagram $D=\{(0,4),(1,3)\}$.  The two tent functions both peak at $t=2$; because $\Lambda_1\geq\Lambda_2$ everywhere, the landscape layers are simply $\lambda_1=\Lambda_1$ and $\lambda_2=\Lambda_2$.  The corresponding scalar statistics read off directly: $A(D)=4$, $L_D=6$, $|I|=2$, and $f_1=f_2=2$, $f_3=f_4=16$.

\begin{center}
\input{figures/section14-landscape-toy}
\captionof{figure}{Landscape construction for $D=\{(0,4),(1,3)\}$.
  Panel~(a): persistence diagram; both birth--death pairs lie above the diagonal.
  Panel~(b): tent functions $\Lambda_1$ (solid) and $\Lambda_2$ (dashed) both peak at $t=2$.
  Because $\Lambda_1\geq\Lambda_2$ everywhere, the first landscape layer $\lambda_1=\Lambda_1$
  (light shading) and second layer $\lambda_2=\Lambda_2$ (darker shading) are the two
  individual tents.}
\label{fig:landscape-toy}
\end{center}

In practice, the persistence diagram fed into this pipeline is derived from the graph-based matrix of \cref{sec:distance-matrix}, not specified by hand.  A periodic signal typically yields one dominant long-lived $H_1$ class: $\lambda_1$ then carries most of the landscape mass, the amplitude $A$ and total persistence $L_D$ are large, and the persistent entropy $E(D)$ is low because a single class dominates the lifetime distribution.  A chaotic signal, by contrast, tends to produce several moderate-persistence classes with more evenly distributed lifetimes, raising $E(D)$ and spreading mass across more landscape layers.  The stability results of \cref{thm:stability,thm:landscape-stability} guarantee that these contrasts are robust to small perturbations in the input, which is what makes them reliable discriminative features.

The next section tests these theoretical expectations directly: whether different graph families succeed on the signal types they are meant to encode, whether the more global diffusion distance of \cref{sec:distance-matrix} improves discrimination, and whether the stability discussion is reflected in noise robustness.

\section{Experimental Validation}
\label{sec:experiments}

The preceding sections defined a modular pipeline---graph construction, distance matrix, Vietoris--Rips filtration, and vectorization---and motivated each design choice theoretically.  This section subjects those choices to three empirical tests on standard time series classification benchmarks:
\begin{enumerate}
  \item \textbf{Graph construction comparison} (\cref{sec:exp-main}): Which of the five network constructions produces the most discriminative persistence features, and does the answer depend on the signal type?
  \item \textbf{Distance-matrix ablation} (\cref{sec:exp-distance}): Does the choice of distance metric on the graph affect classification as much as---or more than---the graph construction itself?
  \item \textbf{Noise robustness} (\cref{sec:exp-noise}): Do persistence-based features degrade gracefully under additive noise, as predicted by the Stability Theorem (\cref{thm:stability})?
\end{enumerate}

\subsection{Experimental Setup}
\label{sec:exp-setup}

\paragraph{Datasets.}
We select twelve univariate datasets from the UCR Time Series Archive~\cite{Dau2019}, spanning six application domains: electrocardiography (ECG200, ECGFiveDays), spectroscopy (Coffee, Meat), motion capture and gait analysis (GunPoint, ToeSegmentation1), engine diagnostics (FordA, FordB), energy monitoring (PowerCons), shape outlines (ArrowHead, Plane), and sensor transients (Trace).  The suite includes both binary and multiclass problems (up to seven classes), with series lengths ranging from $L=96$ to $L=500$ and numbers of series ranging from $N_{\mathrm{obs}}=56$ to $N_{\mathrm{obs}}=4{,}921$.  All datasets are balanced or near-balanced (class imbalance ratio $\leq 1.99$).  \Cref{tab:datasets} summarizes the key properties.
Five additional datasets (Beef, ShapeletSim, SyntheticControl, TwoPatterns, UWaveGestureLibraryX) were considered but excluded because no graph-based pipeline achieved $F_1>0.65$ on any of them.  In each case, the mismatch can be traced to either short series that limit the richness of the reconstructed graph or to class boundaries defined by localized shapelets rather than by global dynamics---regimes in which global topological summaries are inherently uninformative.

For each dataset, we merge the published train and test splits and apply stratified five-fold cross-validation ($K=5$, random seed 42), ensuring that every sample contributes to exactly one held-out test evaluation while each fold trains on a representative $80\%$ of the data.

\begin{table}[t]
\centering
\caption{Summary of the twelve benchmark datasets.  $N_{\mathrm{obs}}$: total number of series (train+test merged); $L$: series length; \#cls: number of classes; imbalance: ratio of the largest to smallest class.}
\label{tab:datasets}
\input{experiments/results/table_datasets}
\end{table}

\paragraph{Pipelines.}
We evaluate five graph constructions defined in \cref{sec:networks}: HVG and NVG (both parameter-free), OPN, CGSSN, and $k$-NN.  Embedding parameters for OPN are selected via MsPE, searching $\tau\in[1,7]$ (a narrower window than the $\tau\in[1,50]$ practical range of \cref{rem:parameter-ranges}, to keep the delay selection numerically stable on the shortest series in the suite) and $n\in[3,7]$; for CGSSN and $k$-NN, $\tau$ is chosen by MI and $n$ by FNN, with a state-space cap of $b^n\leq4{,}096$ for CGSSN (using $b=8$ bins) and $k=5$ for $k$-NN.  Parameters are selected once per dataset (median over a 30-series subset) and shared across all series, ensuring comparable feature scales for the classifier.  \Cref{tab:params} reports the selected values.  All parameter-selection heuristics are invoked through teaspoon's \texttt{parameter\_selection} module; for the OPN we use \texttt{MsPE\_tau} and \texttt{MsPE\_n}, and for the other two we use \texttt{MI\_for\_delay} together with \texttt{FNN\_n}.

In the main experiment (\cref{sec:exp-main}), unweighted graphs (HVG, NVG, $k$-NN) use the shortest-path distance, and weighted graphs (OPN, CGSSN) use lazy-random-walk diffusion distance.  For the latter we set the walk length to $t = \min\!\bigl(\lceil\log_2 |V|\rceil, 10\bigr)$, a size-adaptive and bounded choice that stays consistent with the diameter-based guideline of \cref{sec:distance-matrix}; the shorter walk prevents the over-mixing on small graphs that a full $2\,\mathrm{diam}(G)$ walk can cause, which would otherwise drive the random walk close to its stationary distribution and wash out class-discriminative structure.  The distance-matrix ablation (\cref{sec:exp-distance}) systematically varies this choice.

All four distance matrices are computed through teaspoon's \texttt{network\_tools} module; for reproducibility, we note that the hop-count variant of \cref{sec:distance-matrix} corresponds to the option \texttt{shortest\_weighted\_path}, whose returned value is a hop count rather than a summed path weight despite the name.

\begin{table}[t]
\centering
\caption{Selected embedding parameters ($\tau$, $n$) per dataset.  HVG and NVG are parameter-free.  OPN parameters are chosen via MsPE; CGSSN and $k$-NN share the same MI/$\mathrm{FNN}$-selected values and are therefore reported together.}
\label{tab:params}
\input{experiments/results/table_params}
\end{table}

\paragraph{Vectorization.}
Each persistence diagram (dimensions $H_0$ and $H_1$) is vectorized following the recipe of \cref{sec:feature-assembly}: the mean of the first $J=3$ persistence-landscape layers on a $M=200$-point grid, plus nine scalar summary statistics (persistent entropy, amplitude, total persistence, cardinality, four Adcock--Carlsson coordinates, and $\|\lambda_1\|_1$), yielding a $418$-dimensional feature vector per series.

\paragraph{Classifier.}
We use CatBoost~\cite{Prokhorenkova2018} (400 iterations, depth~5, learning rate 0.05, Logloss for binary and MultiClass for multiclass problems) as the sole classifier to isolate the effect of graph construction and distance choice from classifier complexity.  Zero-variance features are automatically dropped; when all features are constant (e.g., NVG with shortest-path distance on certain datasets), the pipeline falls back to a majority-class baseline.

\paragraph{Metrics.}
We report the macro-averaged $F_1$ score and ROC-AUC (one-vs-rest for multiclass) averaged over five folds.  Bold entries in result tables mark the best pipeline per dataset.

\paragraph{Computational environment.}
All experiments were executed on the TRUBA HPC centre (ARF cluster), using the \texttt{orfoz} queue equipped with dual Intel Xeon Platinum 8480+ processors (112 cores per node) and 256\,GB RAM, running Rocky Linux~9.2.  Feature extraction was parallelized across all available cores using \texttt{joblib}.  \Cref{tab:timing} reports the wall-clock extraction times per dataset and pipeline.

\begin{table}[t]
\centering
\caption{Feature extraction time (seconds) per dataset and pipeline on a 112-core node.}
\label{tab:timing}
\input{experiments/results/table_extraction_time}
\end{table}

\subsection{Graph Construction Comparison}
\label{sec:exp-main}

\Cref{tab:f1,tab:auc} report the classification results per dataset.  Each graph construction is best on a distinct subset of the suite, and the subsets align with the signal structure that the construction was designed to capture.

\emph{CGSSN} achieves the highest $F_1$ on five of twelve datasets: ECG200 (0.618), ECGFiveDays (0.846), PowerCons (0.875), Trace (0.979), and Plane (0.845).  On the first four, the class-discriminative structure resides in the amplitude profile of the signal, which the coarse-grained binning of CGSSN (\cref{def:cgssn}) partitions directly: ECG recordings differ in R-peak height and ST-segment elevation, PowerCons contrasts on--off appliance signatures with gradually ramping loads, and each Trace transient has a distinct amplitude envelope.  Plane is a borderline case, a multiclass shape-outline task on which CGSSN (0.845) and $k$-NN (0.843) are essentially tied; we read this as a boundary between the amplitude-transition and attractor-geometry regimes rather than a clean CGSSN win.

The \emph{OPN} wins on four datasets: Coffee (0.928), FordB (0.770), Meat (0.737), and ToeSegmentation1 (0.739).  In the spectroscopic datasets (Coffee, Meat), the class-discriminative information lies in the ordering of absorption peaks rather than in their absolute intensity; this is precisely the signal that ordinal permutation patterns encode (\cref{def:opn}).  In gait data (ToeSegmentation1), the temporal ordering of heel-to-toe pressure transitions distinguishes classes.  The OPN's deliberate discarding of amplitude information also renders it robust to intensity-calibration differences across samples.

The \emph{$k$-NN} graph achieves the best $F_1$ on ArrowHead (0.692) and GunPoint (0.769).  Both are shape- or motion-capture tasks in which the reconstructed phase-space trajectory encodes the class: gun-drawing produces an asymmetric arc whereas pointing produces a symmetric trajectory.  The $k$-NN graph connects points that are nearby in the Takens embedding (\cref{def:knn}), so its topology reflects attractor geometry rather than temporal succession or amplitude.

The \emph{HVG} achieves $F_1=0.992$ on FordA, a binary classification of normal versus faulty engine sensor signals.  Being parameter-free and amplitude-invariant (\cref{rem:visibility-limitation}), HVG captures the peak/visibility pattern of such signals without being affected by sensor-calibration differences, an advantage in industrial monitoring.  The companion dataset FordB, collected under different operating conditions, tells the opposite story: on FordB, HVG collapses to $F_1=0.337$ while OPN attains the highest $F_1=0.770$.  The two datasets are therefore a useful caution---method--data alignment can be fragile to changes in acquisition conditions even within a nominally similar domain.  Across the remaining datasets, HVG underperforms substantially and often collapses to the majority class ($F_1\approx0.33$), suggesting that its local, horizon-based connectivity is insufficient when class boundaries depend on transitions, ordinal patterns, or global attractor shape.

The \emph{NVG} never outperforms HVG on any dataset under the shortest-path distance, despite its denser connectivity ($O(L\log L)$ edges versus HVG's $O(L)$), and its one-vs-rest ROC-AUC equals $0.50$ exactly on every dataset in the suite---the signature of a persistence diagram whose features are indistinguishable across classes.  The denser connectivity compresses shortest-path distances into a narrow range and yields persistence diagrams with limited inter-class variance.  \Cref{sec:exp-distance} shows that this limitation is almost entirely overcome by switching to diffusion distance.

\begin{table}[t]
\centering
\caption{Mean $F_1$ (macro) $\pm$ std over five folds.  Bold marks the best pipeline per dataset.  The bottom row counts the number of datasets each pipeline wins.}
\label{tab:f1}
\input{experiments/results/table_f1}
\end{table}

\begin{table}[t]
\centering
\caption{Mean ROC-AUC $\pm$ std over five folds (one-vs-rest for multiclass).  Bold marks the best pipeline per dataset.}
\label{tab:auc}
\input{experiments/results/table_auc}
\end{table}

The preceding winners are reported under a single default construction of $\mathbf{D}$ per graph family; the next subsection asks how much of the comparison is driven by that choice.

\subsection{Distance-Matrix Ablation}
\label{sec:exp-distance}

\begin{figure}[t]
\centering
\includegraphics[width=0.72\textwidth]{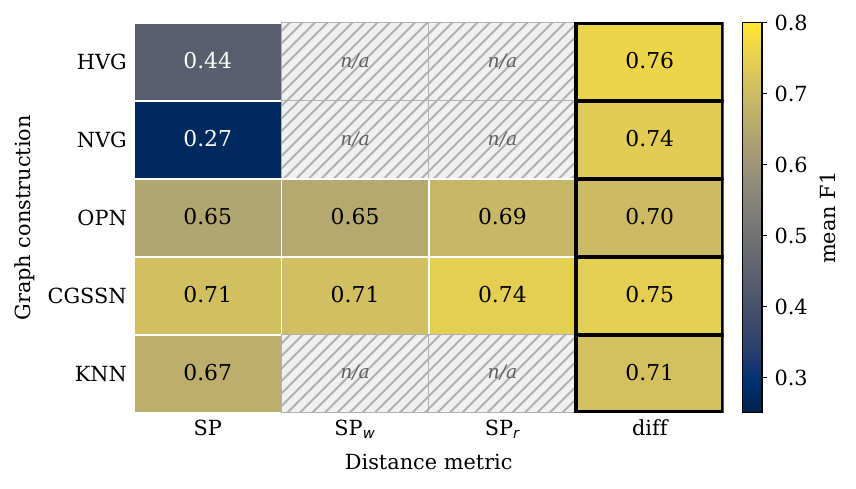}
\caption{Distance-matrix ablation: mean $F_1$ averaged over the twelve datasets, for every (graph, distance) combination.  Rows are graph constructions; columns are distance metrics.  Colour encodes $F_1$ (cividis; darker is lower).  Each row's winning distance is outlined in a bold black border.  Hatched cells (\emph{n/a}) denote combinations that are not applicable---unweighted graphs (HVG, NVG, $k$-NN) do not admit weight-based distances.  Diffusion distance is the row-winner for every graph family, with the largest gains over shortest-path distance on the visibility graphs.  Distance labels: SP, shortest unweighted path; SP$_w$, hop count on a reciprocal-weight-optimal path; SP$_r$, reciprocal-weight shortest path; diff, lazy-random-walk diffusion distance.}
\label{fig:distance-heatmap}
\end{figure}

\Cref{sec:distance-matrix} defined four ways to construct the pairwise matrix $\mathbf{D}$ from a graph: shortest unweighted path, hop count on a reciprocal-weight-optimal path, reciprocal-weight shortest path, and diffusion distance.  The main experiment intentionally fixed one default choice of $\mathbf{D}$ per graph family.  Here we hold the graph construction fixed and vary the construction of $\mathbf{D}$, testing all applicable combinations (two for unweighted graphs, four for weighted graphs; fourteen combinations in total).  \Cref{fig:distance-heatmap} summarises the result as a heatmap of mean $F_1$ averaged across the twelve datasets.

For every graph construction, \emph{diffusion distance is the row-wise winner} in \cref{fig:distance-heatmap}.  Replacing shortest-path distance by diffusion distance yields a modest improvement on the weighted graphs---$5.1\%$ for CGSSN ($0.710\to0.746$), $6.8\%$ for $k$-NN ($0.665\to0.710$), and $8.1\%$ for OPN ($0.645\to0.697$)---but a far larger one on the visibility graphs: HVG improves by $74\%$ ($0.437\to0.761$) and NVG by $173\%$ ($0.271\to0.739$).  With diffusion distance, NVG becomes competitive and even wins two datasets outright (ECGFiveDays and ToeSegmentation1), recovering from its collapse under shortest-path distance.

A single-dataset illustration sharpens the point.  On ECGFiveDays (\cref{tab:distance-ecgfivedays}), the best graph under shortest-path distance is HVG at $F_1=0.627$, yet switching the distance on the same NVG graph---from shortest-path to diffusion---raises $F_1$ to $0.930$.  The gain attributable to changing the distance exceeds any gain attainable by changing the graph at a fixed distance, so the construction of $\mathbf{D}$ is a first-order design decision rather than a technical detail.

\begin{table}[t]
\centering\small
\caption{ECGFiveDays as a single-dataset illustration of the claim in \cref{sec:exp-distance}: switching the distance metric on a given graph yields a larger $F_1$ change than switching the graph at a fixed distance.  Best value in bold.}
\label{tab:distance-ecgfivedays}
\begin{tabular}{lcc}
\toprule
Graph & $F_1$ (shortest path) & $F_1$ (diffusion) \\
\midrule
HVG   & 0.627 & 0.884 \\
NVG   & 0.332 & \textbf{0.930} \\
OPN   & 0.498 & 0.686 \\
CGSSN & 0.801 & 0.846 \\
$k$-NN & 0.568 & 0.748 \\
\bottomrule
\end{tabular}
\end{table}

\subsection{Noise Robustness}
\label{sec:exp-noise}

\Cref{thm:stability}, due to Cohen-Steiner, Edelsbrunner, and Harer~\cite{CohenSteiner2007}, guarantees that persistence diagrams are Lipschitz-continuous with respect to bounded perturbations of the filtration function.  We do not claim to strengthen that theorem; rather, we examine empirically whether the pipeline exhibits the qualitative behaviour the theorem predicts.  To this end, we inject additive Gaussian noise at five signal-to-noise ratios ($\mathrm{SNR}\in\{\infty, 20, 10, 5, 0\}\,$dB) into two representative datasets---Coffee (where OPN wins on clean data) and PowerCons (where CGSSN wins)---and re-run the full pipeline.  Here the decibel-scale ratio is defined by $\mathrm{SNR}_{\mathrm{dB}} = 10\log_{10}(P_{\mathrm{signal}}/P_{\mathrm{noise}})$, so larger SNR means a cleaner signal, $0$\,dB means equal signal and noise power, and $\infty$ corresponds to the noiseless baseline.  Embedding parameters are selected on the clean signal and held fixed across all noise levels, isolating the effect of input perturbation from parameter instability.

\begin{figure}[t]
\centering
\includegraphics[width=0.85\textwidth]{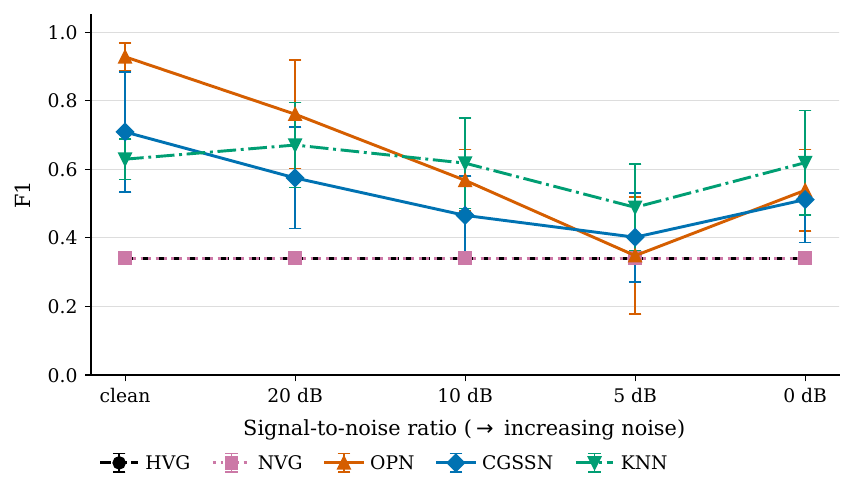}
\caption{Noise robustness on Coffee: $F_1$ versus SNR for each pipeline.  Error bars show the standard deviation across five folds.  OPN loses its clean-data advantage at moderate noise, while CGSSN degrades most gracefully.}
\label{fig:noise-coffee}
\end{figure}

\begin{figure}[t]
\centering
\includegraphics[width=0.85\textwidth]{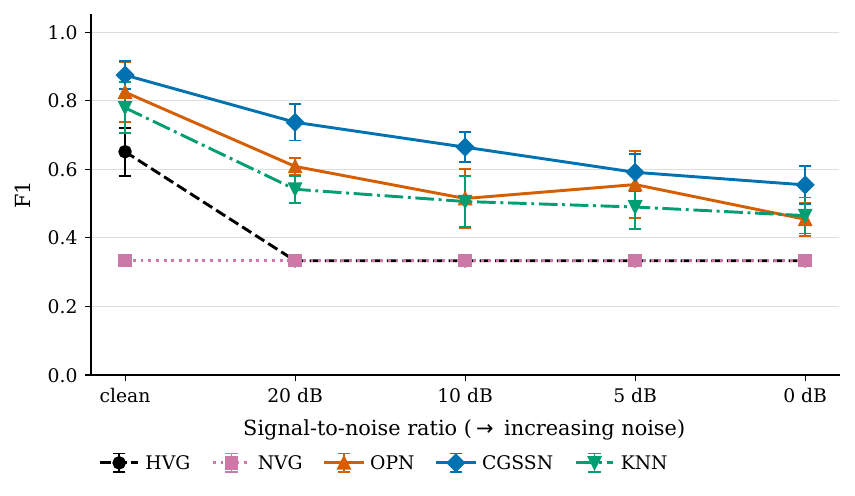}
\caption{Noise robustness on PowerCons: $F_1$ versus SNR.  CGSSN maintains the highest $F_1$ at every noise level.  HVG collapses to the majority-class baseline even on clean data.}
\label{fig:noise-powercons}
\end{figure}

\Cref{fig:noise-coffee,fig:noise-powercons} plot $F_1$ against SNR for each pipeline, and three patterns emerge.

First, all competitive pipelines \emph{degrade gracefully}: the curves slope downward smoothly as noise increases, without an abrupt performance cliff.  This is the qualitative behaviour predicted by the Lipschitz bound of \cref{thm:stability}.

Second, the \emph{CGSSN is the most noise-robust} among the evaluated constructions.  On Coffee, CGSSN retains $72\%$ of its clean $F_1$ at $0$\,dB SNR (from $0.709$ to $0.512$); on PowerCons, the corresponding retention is $63\%$ (from $0.875$ to $0.554$).  In both datasets, CGSSN maintains the highest absolute $F_1$ among the evaluated pipelines at every noise level below $20$\,dB.  The amplitude binning of the CGSSN (\cref{def:cgssn}) provides a natural quantisation that absorbs small additive perturbations without altering the coarse transition structure.

Third, the \emph{OPN is the most noise-sensitive} competitive method.  On Coffee, its $F_1$ falls from $0.928$ to $0.539$ at $0$\,dB (a retention of $58\%$), yielding its clean-data advantage to CGSSN.  This behaviour follows from the construction itself: additive noise can reorder nearby values, which in turn scrambles the ordinal permutation patterns on which the OPN relies.  The CGSSN should therefore be preferred over the OPN when the noise level is appreciable, even on tasks where the OPN is superior in the noiseless regime.

\subsection{Discussion}
\label{sec:exp-discussion}

Taken together, the three experiments yield three complementary lessons.  First, no single graph construction dominates universally; the optimal choice depends on where the class-discriminative information resides in the signal, and the method-data alignments reported in \cref{sec:exp-main} can serve as a principled starting point in new applications.  Second, the construction of $\mathbf{D}$ is a first-order design choice rather than a technical afterthought: the empirical advantage of diffusion distance is consistent with the more global random-walk geometry described in \cref{sec:distance-matrix}, especially for dense visibility graphs whose shortest-path distances collapse into a narrow range.  Third, the Stability Theorem of \cref{sec:stability} is reflected empirically, though not yet proved end to end: the noise experiment shows gradual rather than catastrophic loss of accuracy, which is the qualitative behaviour one would expect when persistence diagrams vary Lipschitz-continuously under perturbation.

These observations set the stage for the concluding remarks on practical use and remaining open problems.

\section{Conclusion and Future Directions}
\label{sec:conclusion}

The framework developed in this study provides a unified route from univariate time series to persistence-based features through complex networks.  By separating graph construction, graph-derived dissimilarity, persistence computation, and vectorization, the framework makes it possible to compare competing representations on common downstream ground rather than mixing changes in the graph with changes in the persistence or learning stages.

The experiments suggest that no single graph construction is uniformly best; what matters is the relation between the signal and the aspect of structure emphasized by the graph.  Within the benchmark suite studied here, CGSSN combined with diffusion distance is the strongest overall default, particularly on ECG, energy, and transient data.  OPN remains attractive when class information is carried mainly by ordinal structure rather than amplitude, while $k$-NN graphs are preferable when the reconstructed geometry itself is discriminative.  Visibility graphs are less reliable as default classifiers, but they retain practical value when parameter-free construction and direct interpretability from the signal are primary considerations.  For very long series, the symbolic OPN and CGSSN are also plausible starting points on complexity grounds because they keep the vertex set manageable.  Across the graph types tested here, the choice of graph-derived distance is not secondary: diffusion distance consistently improves on shortest-path alternatives and is therefore the preferable starting point whenever the additional computation is acceptable.

Several limitations remain.  The transition and proximity constructions depend on embedding and graph parameters whose heuristics do not always agree across signal classes, and the computational burden of large visibility or $k$-NN graphs remains substantial despite the compression offered by symbolic methods.  The study is restricted to univariate series, so extension to multivariate settings requires further work on both network design and topological interpretation.  On the theoretical side, the Stability Theorem (\cref{thm:stability}) controls the persistence stage once the filtration is specified, but the stability of the full chain from raw signal to network to graph-derived matrix is still open.  Finally, vectorization necessarily discards some diagram-level information, so the interaction between network-based diagrams and specific vectorization schemes deserves more systematic study.

Taken together, these results suggest that network-based persistent homology is best understood as a modular representation framework rather than as a universal classifier.  Its practical strength lies in allowing different graph constructions to encode different kinds of temporal structure while retaining a common topological downstream analysis.

\section*{Acknowledgments}
\label{sec:acknowledgments}

The author is deeply grateful to his wife, Berna, and his children, Masal and Kuzey, for their patience and understanding during the preparation of this manuscript. 

The numerical calculations reported in this study were fully performed at TÜBİTAK ULAKBİM, High Performance and Grid Computing Center (TRUBA resources).

\begingroup
\let\cleardoublepage\clearpage
\bibliographystyle{plainnat}
\bibliography{references}
\endgroup

\end{document}

%% file: figures/common-route-flow.tex
\begin{tikzpicture}[
  font=\footnotesize,
  >=Latex,
  node distance=4.5mm,
  box/.style={
    draw=black!70,
    rounded corners=2pt,
    fill=black!3,
    align=left,
    text width=0.70\textwidth,
    inner sep=5pt
  },
  downstream/.style={
    box,
    fill=black!7
  },
  route/.style={-Latex, semithick, draw=black!70},
  focus/.style={draw=black!55, dashed, rounded corners=4pt, inner sep=6pt}
]
\node[box] (series) {\textbf{Time-series input.} Univariate signal $x=[x_1,\ldots,x_L]$.};
\node[box, below=of series] (construction) {\textbf{Network construction.} Visibility, transition, or proximity map.};
\node[box, below=of construction] (graph) {\textbf{Graph representation.} $G=(V,E)$ or $G=(V,E,\omega)$.};
\node[downstream, below=8mm of graph] (matrix) {\textbf{Graph-based matrix $\mathbf{D}$.} Shortest path, reciprocal-weight path, diffusion distance, or related separation measure.};
\node[downstream, below=of matrix] (filtration) {\textbf{Flag / Vietoris--Rips filtration.} Build $K_\varepsilon$ by thresholding $\mathbf{D}$ across scales.};
\node[downstream, below=of filtration] (homology) {\textbf{Persistent homology.} Track $H_0$ and $H_1$ across the filtration.};
\node[downstream, below=of homology] (diagram) {\textbf{Persistence diagrams.} Output $\Dgm_0$ and $\Dgm_1$.};

\draw[route] (series) -- (construction);
\draw[route] (construction) -- (graph);
\draw[route] (graph) -- (matrix);
\draw[route] (matrix) -- (filtration);
\draw[route] (filtration) -- (homology);
\draw[route] (homology) -- (diagram);

\node[focus, fit=(matrix)(filtration)(homology)(diagram)] (standardized) {};
\node[fill=white, inner sep=1.5pt, anchor=south west, xshift=4pt, yshift=3pt] at (standardized.north west) {\scriptsize\textbf{Standardized downstream route}};
\end{tikzpicture}

%% file: figures/rips-toy-filtration.tex
\resizebox{0.94\textwidth}{!}{%
\begin{tikzpicture}[x=1cm,y=1cm, >=Latex, baseline]
\tikzset{
  paneltitle/.style={font=\normalfont\footnotesize},
  panelsub/.style={font=\scriptsize},
  topvertex/.style={circle, draw=black, fill=white, minimum size=15pt, inner sep=0pt, font=\scriptsize},
  filtvertex/.style={circle, draw=black, fill=white, minimum size=13pt, inner sep=0pt, font=\scriptsize}
}
\useasboundingbox (0,-0.85) rectangle (13.5,6.6);

\begin{scope}[shift={(1.25,4.75)}]
\node[topvertex] (g1) at (-0.55,0.55) {1};
\node[topvertex] (g2) at (0.55,0.55) {2};
\node[topvertex] (g3) at (0.55,-0.55) {3};
\node[topvertex] (g4) at (-0.55,-0.55) {4};
\draw[semithick] (g1) -- (g2) -- (g3) -- (g4) -- (g1);
\end{scope}
\node[paneltitle] at (1.25,3.18) {(a) Graph};

\node[anchor=north] at (4.75,5.84) {{\scriptsize
  \setlength{\tabcolsep}{4pt}
  \renewcommand{\arraystretch}{1.05}
  \begin{tabular}{c|cccc}
        & 1 & 2 & 3 & 4 \\\hline
    1 & 0 & 1 & 2 & 1 \\
    2 & 1 & 0 & 1 & 2 \\
    3 & 2 & 1 & 0 & 1 \\
    4 & 1 & 2 & 1 & 0
  \end{tabular}
}};
\node[paneltitle] at (4.75,3.18) {(b) Distance Matrix};

\begin{scope}[shift={(7.75,4.7)}]
\node[panelsub] at (0,0.85) {$\varepsilon = 0$};
\node[filtvertex] (f01) at (-0.38,0.25) {1};
\node[filtvertex] (f02) at (0.38,0.25) {2};
\node[filtvertex] (f03) at (0.38,-0.45) {3};
\node[filtvertex] (f04) at (-0.38,-0.45) {4};
\end{scope}

\node[font=\large] at (8.90,4.58) {$\subseteq$};

\begin{scope}[shift={(10.05,4.7)}]
\node[panelsub] at (0,0.85) {$\varepsilon = 1$};
\node[filtvertex] (f11) at (-0.38,0.25) {1};
\node[filtvertex] (f12) at (0.38,0.25) {2};
\node[filtvertex] (f13) at (0.38,-0.45) {3};
\node[filtvertex] (f14) at (-0.38,-0.45) {4};
\draw[semithick] (f11) -- (f12) -- (f13) -- (f14) -- (f11);
\end{scope}

\node[font=\large] at (11.20,4.58) {$\subseteq$};

\begin{scope}[shift={(12.35,4.7)}]
\node[panelsub] at (0,0.85) {$\varepsilon = 2$};
\fill[black!8] (0, -0.1) circle [radius=0.8];
\node[filtvertex] (f21) at (-0.38,0.25) {1};
\node[filtvertex] (f22) at (0.38,0.25) {2};
\node[filtvertex] (f23) at (0.38,-0.45) {3};
\node[filtvertex] (f24) at (-0.38,-0.45) {4};
\draw[semithick] (f21) -- (f22) -- (f23) -- (f24) -- (f21);
\draw[semithick] (f21) -- (f23);
\draw[semithick] (f22) -- (f24);
\end{scope}
\node[paneltitle] at (10.20,3.18) {(c) Filtration};

\node[paneltitle] at (4.55,-0.52) {(d) Barcode};
\node[paneltitle] at (8.95,-0.52) {(e) Diagram};

\begin{scope}[shift={(2.85,0.00)}]
\path[use as bounding box] (0,0) rectangle (3.45,2.45);
\draw[->] (0.55,0.35) -- (3.15,0.35) node[right] {$\varepsilon$};
\foreach \x in {0,1,2} {
  \draw (0.55+\x,0.42) -- (0.55+\x,0.28);
  \node[below, font=\scriptsize] at (0.55+\x,0.28) {\x};
}
\node[anchor=east, font=\small] at (0.35,1.55) {$H_0$};
\draw[gray, line width=0.95pt] (0.55,1.85) -- (1.55,1.85);
\draw[gray, line width=0.95pt] (0.55,1.55) -- (1.55,1.55);
\draw[gray, line width=0.95pt] (0.55,1.25) -- (1.55,1.25);
\draw[black, line width=1.05pt, ->] (0.55,0.95) -- (3.10,0.95);
\node[anchor=east, font=\small] at (0.35,0.55) {$H_1$};
\draw[black, line width=1.05pt] (1.55,0.55) -- (2.55,0.55);
\end{scope}

\begin{scope}[shift={(7.75,0.00)}]
\path[use as bounding box] (0,0) rectangle (2.95,2.45);
\draw[->] (0.35,0.15) -- (2.75,0.15) node[right] {$b$};
\draw[->] (0.35,0.15) -- (0.35,2.45);
\node[anchor=south east, font=\small] at (0.30,2.48) {$d$};
\draw[gray] (0.35,0.15) -- (2.45,2.25);
\foreach \x in {1,2} {
  \draw (0.35+0.85*\x,0.20) -- (0.35+0.85*\x,0.10);
  \node[below, font=\scriptsize] at (0.35+0.85*\x,0.10) {\x};
}
\foreach \y/\label in {1/1,2/2} {
  \draw (0.40,0.15+0.85*\y) -- (0.30,0.15+0.85*\y);
  \node[left, font=\scriptsize] at (0.30,0.15+0.85*\y) {\label};
}
\node[left, font=\scriptsize] at (0.30,2.18) {$\infty$};
\fill[gray] (0.45,0.95) circle (1.1pt);
\fill[gray] (0.55,1.05) circle (1.1pt);
\fill[gray] (0.62,0.88) circle (1.1pt);
\fill[gray] (0.50,1.95) circle (1.1pt);
\draw[gray, ->] (0.50,1.95) -- (0.50,2.22);
\fill[black] (1.20,1.85) circle (1.2pt);
\end{scope}
\end{tikzpicture}%
}

%% file: figures/section13-visibility-toy.tex
%
%
%
%
\resizebox{0.70\textwidth}{!}{%
\begin{tikzpicture}[x=1cm,y=1cm, >=Latex, baseline]
\tikzset{
  paneltitle/.style={font=\normalfont\large, anchor=north},
  tspt/.style={circle, fill=black, inner sep=1.5pt},
  gnode/.style={circle, draw=black, fill=white, minimum size=18pt, inner sep=0pt,
                font=\small, line width=0.70pt},
  gedge/.style={line width=0.70pt, shorten >=3pt, shorten <=3pt},
  missingedge/.style={line width=1.50pt, dash pattern=on 5pt off 3pt,
                      draw=black!55, shorten >=3pt, shorten <=3pt}
}
\useasboundingbox (-0.20,-0.60) rectangle (15.40,11.20);


\begin{scope}[shift={(1.20,6.60)}]
  \draw[->](-0.15,0)--(4.70,0)node[right,font=\small]{$t$};
  \draw[->](0,-0.15)--(0,2.90)node[above,font=\small]{$x_t$};
  \foreach \i/\xi in {1/0.82,2/1.64,3/2.46,4/3.28,5/4.10}{
    \draw(\xi,0.06)--(\xi,-0.06)node[below=2pt,font=\small]{\i};
  }
  \foreach \y/\lbl in {0.50/1,1.00/2,1.50/3,2.00/4,2.50/5}{
    \draw(0.06,\y)--(-0.06,\y)node[left=2pt,font=\small]{\lbl};
  }
  \draw[line width=0.70pt](0.82,0.50)--(1.64,1.00)--(2.46,0.50)--(3.28,2.50)--(4.10,1.00);
  \node[tspt]at(0.82,0.50){};
  \node[tspt]at(1.64,1.00){};
  \node[tspt]at(2.46,0.50){};
  \node[tspt]at(3.28,2.50){};
  \node[tspt]at(4.10,1.00){};
\end{scope}
\node[paneltitle]at(3.80,5.32){(a)~Signal $x=[1,2,1,5,2]$};

\begin{scope}[shift={(8.80,7.90)}]
  \foreach \i/\xi in {1/0.0,2/1.3,3/2.6,4/3.9,5/5.2}{
    \node[gnode](n\i)at(\xi,0){\i};
  }
  \foreach \i/\j in {1/2,2/3,3/4,4/5}{\draw[gedge](n\i)--(n\j);}
  \draw[gedge](n2)..controls+(90:0.80)and+(90:0.80)..(n4);
  \draw[gedge](n1)..controls+(90:1.55)and+(90:1.55)..(n4);
\end{scope}
\node[paneltitle]at(11.60,5.32){(b)~NVG (linear layout)};


\begin{scope}[shift={(0.70,0.20)}]
  \node[gnode](o4)at(2.50,4.10){4};
  \node[gnode](o3)at(0.79,2.86){3};
  \node[gnode](o2)at(1.44,0.84){2};
  \node[gnode](o1)at(3.56,0.84){1};
  \node[gnode](o5)at(4.21,2.86){5};
  \draw[gedge](o1)--(o2);
  \draw[gedge](o2)--(o3);
  \draw[gedge](o3)--(o4);
  \draw[gedge](o4)--(o5);
  \draw[gedge](o2)--(o4);   
  \draw[gedge](o1)--(o4);   
\end{scope}
\node[paneltitle]at(3.50,-0.18){(c)~NVG (graph layout)};

\begin{scope}[shift={(8.60,0.20)}]
  \node[gnode](h4)at(2.50,4.10){4};
  \node[gnode](h3)at(0.79,2.86){3};
  \node[gnode](h2)at(1.44,0.84){2};
  \node[gnode](h1)at(3.56,0.84){1};
  \node[gnode](h5)at(4.21,2.86){5};
  \draw[gedge](h1)--(h2);
  \draw[gedge](h2)--(h3);
  \draw[gedge](h3)--(h4);
  \draw[gedge](h4)--(h5);
  \draw[gedge](h2)--(h4);              
  \draw[missingedge](h1)--(h4);        
\end{scope}
\node[paneltitle]at(11.30,-0.18){(d)~HVG (graph layout)};

\end{tikzpicture}%
}

%% file: figures/section13-opn-toy.tex
%
%
%
%
\resizebox{0.88\textwidth}{!}{%
\begin{tikzpicture}[x=1cm,y=1cm, >=Latex, baseline]
\tikzset{
  paneltitle/.style={font=\normalfont\small, anchor=north},
  tspt/.style={circle, fill=black, inner sep=1.5pt},
  opnnode/.style={draw=black, fill=white, circle,
                  minimum size=0.56cm, inner sep=1pt,
                  font=\small},
  gedge/.style={line width=0.70pt, shorten >=4pt, shorten <=4pt},
  heavyedge/.style={line width=1.80pt, shorten >=4pt, shorten <=4pt}
}
\useasboundingbox (-0.20,-0.70) rectangle (15.00,5.10);

\begin{scope}[shift={(0.50,0.90)}]
  \fill[gray!12] (0.35,-0.12) rectangle (2.45,3.10);
  \draw[->](-0.15,0)--(7.00,0)node[right,font=\small]{$t$};
  \draw[->](0,-0.15)--(0,3.40)node[above,font=\small]{$x_t$};
  \foreach \i/\xi in {1/0.70,2/1.40,3/2.10,4/2.80,5/3.50,6/4.20,7/4.90,8/5.60,9/6.30}{
    \draw(\xi,0.06)--(\xi,-0.06)node[below=2pt,font=\small]{\i};
  }
  \foreach \y/\lbl in {0.55/1,1.10/2,1.65/3,2.20/4,2.75/5}{
    \draw(0.06,\y)--(-0.06,\y)node[left=2pt,font=\small]{\lbl};
  }
  \draw[line width=0.70pt]
    (0.70,0.55)--(1.40,1.65)--(2.10,1.10)--(2.80,2.75)
    --(3.50,2.20)--(4.20,1.10)--(4.90,2.20)--(5.60,0.55)--(6.30,1.65);
  \foreach \xi/\yi in {
      0.70/0.55,1.40/1.65,2.10/1.10,2.80/2.75,
      3.50/2.20,4.20/1.10,4.90/2.20,5.60/0.55,6.30/1.65}{
    \node[tspt] at (\xi,\yi) {};
  }
  \draw[decorate, decoration={brace, amplitude=7pt, raise=7pt}]
    (2.10,-0.28) -- (0.70,-0.28)
    node[midway, below=9pt, font=\scriptsize]{$n{=}3$};
\end{scope}
\node[paneltitle]at(3.90,0.06){(a)~Signal $x=[1,3,2,5,4,2,4,1,3]$};

\begin{scope}[shift={(8.80,1.40)}]
  \node[opnnode](pa) at (0.00, 0.00) {$\pi_1$};
  \node[opnnode](pb) at (1.80, 0.00) {$\pi_2$};
  \node[opnnode](pd) at (3.60, 0.00) {$\pi_4$};
  \node[opnnode](pc) at (5.40, 0.00) {$\pi_3$};
  \node[opnnode](pe) at (1.80, 1.80) {$\pi_5$};
  \draw[heavyedge](pa) -- node[above=1pt, font=\scriptsize]{$2$} (pb);
  \draw[gedge](pa) -- node[left=1pt, font=\scriptsize]{$1$} (pe);
  \draw[gedge](pb) -- node[right=1pt, font=\scriptsize]{$1$} (pe);
  \draw[gedge](pb) -- node[above=1pt, font=\scriptsize]{$1$} (pd);
  \draw[gedge](pd) -- node[above=1pt, font=\scriptsize]{$1$} (pc);
\end{scope}
\node[paneltitle]at(11.50,0.06){(b)~OPN ($n{=}3$,~$\tau{=}1$)};

\end{tikzpicture}%
}

%% file: figures/section13-cgssn-toy.tex
%
%
%
%
%
\resizebox{0.88\textwidth}{!}{%
\begin{tikzpicture}[x=1cm,y=1cm, >=Latex, baseline]
\tikzset{
  paneltitle/.style={font=\normalfont\small, anchor=north},
  tspt/.style={circle, fill=black, inner sep=1.5pt},
  cgssnnode/.style={draw=black, fill=white, circle,
                    minimum size=0.56cm, inner sep=1pt, font=\small},
  gedge/.style={line width=0.70pt, shorten >=4pt, shorten <=4pt},
  heavyedge/.style={line width=1.80pt, shorten >=4pt, shorten <=4pt}
}
\useasboundingbox (-0.20,-0.70) rectangle (15.00,5.10);

\begin{scope}[shift={(0.60,0.90)}]
  \fill[gray!10] (0.0,0.0) rectangle (7.20,1.80);
  \fill[gray!22] (0.0,1.80) rectangle (7.20,3.50);
  \draw[dashed, gray!70, line width=0.60pt] (0.0,1.80) -- (7.20,1.80);
  \node[font=\scriptsize, gray!60, anchor=east] at (7.20,0.90) {bin~0};
  \node[font=\scriptsize, gray!60, anchor=east] at (7.20,2.60) {bin~1};
  \draw[->](-0.15,0)--(7.00,0)node[right,font=\small]{$t$};
  \draw[->](0,-0.15)--(0,3.60)node[above,font=\small]{$x_t$};
  \foreach \i/\xi in {1/0.70,2/1.40,3/2.10,4/2.80,5/3.50,6/4.20,7/4.90,8/5.60,9/6.30}{
    \draw(\xi,0.06)--(\xi,-0.06)node[below=2pt,font=\small]{\i};
  }
  \foreach \y/\lbl in {0.40/1,1.60/4,2.00/5,3.20/8}{
    \draw(0.06,\y)--(-0.06,\y)node[left=2pt,font=\small]{\lbl};
  }
  \node[right=2pt, font=\scriptsize, gray!70] at (7.20,1.80) {$4.5$};
  \draw[line width=0.70pt]
    (0.70,0.40)--(1.40,1.60)--(2.10,2.80)--(2.80,0.80)
    --(3.50,2.00)--(4.20,3.20)--(4.90,1.20)--(5.60,2.40)--(6.30,0.40);
  \foreach \xi/\yi in {
    0.70/0.40,1.40/1.60,2.10/2.80,2.80/0.80,
    3.50/2.00,4.20/3.20,4.90/1.20,5.60/2.40,6.30/0.40}{
    \node[tspt] at (\xi,\yi) {};
  }
  \draw[decorate, decoration={brace, amplitude=7pt, raise=7pt}]
    (2.10,-0.28) -- (0.70,-0.28)
    node[midway, below=9pt, font=\scriptsize]{$n{=}3$};
\end{scope}
\node[paneltitle]at(4.00,0.06){(a)~Signal $x=[1,4,7,2,5,8,3,6,1]$};

\begin{scope}[shift={(8.70,1.10)}]
  \node[cgssnnode](s5) at (0.00, 0.00) {$5$};
  \node[cgssnnode](s3) at (1.80, 0.00) {$3$};
  \node[cgssnnode](s6) at (3.60, 0.00) {$6$};
  \node[cgssnnode](s4) at (5.40, 0.00) {$4$};
  \node[cgssnnode](s7) at (3.60, 1.80) {$7$};
  \draw[heavyedge](s3) -- node[above=1pt, font=\scriptsize]{$2$} (s6);
  \draw[gedge](s5) -- node[above=1pt, font=\scriptsize]{$1$} (s3);
  \draw[gedge](s6) -- node[right=1pt, font=\scriptsize]{$1$} (s7);
  \draw[gedge](s6) -- node[above=1pt, font=\scriptsize]{$1$} (s4);
  \draw[gedge](s4) -- node[right=1pt, font=\scriptsize]{$1$} (s7);
\end{scope}
\node[paneltitle]at(12.00,0.06){(b)~CGSSN ($n{=}3$,~$\tau{=}1$,~$b{=}2$)};

\end{tikzpicture}%
}

%% file: figures/section13-knn-toy.tex
%
%
\resizebox{0.88\textwidth}{!}{%
\begin{tikzpicture}[
  gedge/.style={draw=black!70, line width=1.00pt},
  tspt/.style={circle, fill=black, inner sep=0pt, minimum size=3.5pt},
  vnode/.style={circle, draw=black!80, fill=white, minimum size=0.60cm,
                inner sep=0pt, font=\small},
]

\begin{scope}[shift={(1.00,0.65)}]

  \fill[gray!12] (0.5*0.68, -0.12) rectangle (3.5*0.68, 5*0.52+0.15);

  \draw[->] (-0.05,0) -- (7*0.68+0.38,0) node[right, font=\small] {$t$};
  \draw[->] (0,-0.05) -- (0, 5*0.52+0.40) node[above, font=\small] {$x_t$};

  \foreach \yv in {1,...,5} {
    \draw (0,\yv*0.52) -- (-0.06,\yv*0.52);
    \node[left, font=\scriptsize] at (-0.07,\yv*0.52) {\yv};
  }

  \foreach \ti in {1,...,7} {
    \draw (\ti*0.68, 0) -- (\ti*0.68, -0.06);
    \node[below, font=\scriptsize] at (\ti*0.68, -0.05) {\ti};
  }

  \draw[line width=0.70pt]
    (1*0.68, 1*0.52) --
    (2*0.68, 2*0.52) --
    (3*0.68, 1*0.52) --
    (4*0.68, 5*0.52) --
    (5*0.68, 2*0.52) --
    (6*0.68, 3*0.52) --
    (7*0.68, 1*0.52);

  \foreach \ti/\yv in {1/1, 2/2, 3/1, 4/5, 5/2, 6/3, 7/1}
    \node[tspt] at (\ti*0.68, \yv*0.52) {};

  \draw[decorate, decoration={brace, amplitude=7pt, raise=7pt, mirror}]
    (1*0.68, -0.30) -- (3*0.68, -0.30)
    node[midway, below=9pt, font=\scriptsize] {$n{=}3$};

  \node[font=\normalfont\small, anchor=north] at (4*0.68, -1.05)
    {(a) Signal $x = [1,2,1,5,2,3,1]$};

\end{scope}

\begin{scope}[shift={(7.80,0.65)}]

  \draw[->] (-0.10,0) -- (6*0.52+0.38,0) node[right, font=\small] {$\xi_1$};
  \draw[->] (0,-0.10) -- (0, 5*0.52+0.40) node[above, font=\small] {$\xi_2$};

  \foreach \v in {1,...,5} {
    \draw (\v*0.52, 0)    -- (\v*0.52, -0.06);
    \node[below, font=\scriptsize] at (\v*0.52, -0.05) {\v};
    \draw (0, \v*0.52)    -- (-0.06, \v*0.52);
    \node[left, font=\scriptsize]  at (-0.07, \v*0.52) {\v};
  }

  \coordinate (v1) at (1*0.52, 2*0.52);  
  \coordinate (v2) at (2*0.52, 1*0.52);  
  \coordinate (v3) at (1*0.52, 5*0.52);  
  \coordinate (v4) at (5*0.52, 2*0.52);  
  \coordinate (v5) at (2*0.52, 3*0.52);  

  \draw[gedge] (v1) -- (v2);
  \draw[gedge] (v1) -- (v3);
  \draw[gedge] (v1) -- (v5);
  \draw[gedge] (v2) -- (v4);
  \draw[gedge] (v3) -- (v5);
  \draw[gedge] (v4) -- (v5);

  \node[vnode] at (v1) {$v_1$};
  \node[vnode] at (v2) {$v_2$};
  \node[vnode] at (v3) {$v_3$};
  \node[vnode] at (v4) {$v_4$};
  \node[vnode] at (v5) {$v_5$};

  \node[font=\normalfont\small, anchor=north] at (1.75, -1.05)
    {(b) $2$-NN graph ($n{=}3$, $\tau{=}1$)};

\end{scope}

\useasboundingbox (-0.20,-1.40) rectangle (15.00,4.40);

\end{tikzpicture}%
}

%% file: figures/section14-landscape-toy.tex
\resizebox{0.88\textwidth}{!}{%
\begin{tikzpicture}[x=1cm,y=1cm, >=Latex, baseline]
\tikzset{
  paneltitle/.style={font=\normalfont\small},
  axislabel/.style={font=\scriptsize},
  ptlabel/.style={font=\scriptsize, inner sep=1.5pt}
}
\useasboundingbox (-0.20,-1.00) rectangle (15.00,4.20);

\begin{scope}[shift={(1.20, 0.60)}, x=0.65cm, y=0.65cm]
  \draw[->] (-0.10, 0) -- (4.50, 0) node[right, axislabel] {$b$};
  \draw[->] (0, -0.10) -- (0, 4.50) node[above, axislabel] {$d$};
  \draw[gray, thin] (0,0) -- (4,4);
  \foreach \x in {1,2,3,4} {
    \draw (\x,0.06) -- (\x,-0.06);
    \node[below, axislabel] at (\x, -0.06) {\x};
  }
  \foreach \y in {1,2,3,4} {
    \draw (0.06,\y) -- (-0.06,\y);
    \node[left, axislabel] at (-0.06,\y) {\y};
  }
  \fill[black] (0,4) circle [radius=3pt];
  \node[ptlabel, right=2pt] at (0,4) {$(b_1,d_1){=}(0,4)$};
  \fill[black!50] (1,3) circle [radius=3pt];
  \node[ptlabel, above right] at (1,3) {$(b_2,d_2){=}(1,3)$};
  \draw[gray, dashed, very thin] (0,4) -- (2,2);
  \draw[gray, dashed, very thin] (1,3) -- (2,2);
\end{scope}
\node[paneltitle] at (2.65, -0.60) {(a) Persistence diagram};

\begin{scope}[shift={(7.80, 0.60)}]
  \draw[->] (-0.10, 0) -- (5.10, 0) node[right, axislabel] {$t$};
  \draw[->] (0, -0.10) -- (0, 2.70) node[above, axislabel] {};
  \foreach \x in {1,2,3,4} {
    \draw (\x,0.06) -- (\x,-0.06);
    \node[below, axislabel] at (\x,-0.06) {\x};
  }
  \foreach \y in {1,2} {
    \draw (0.06,\y) -- (-0.06,\y);
    \node[left, axislabel] at (-0.06,\y) {\y};
  }

  \draw[black, semithick] (0,0) -- (2,2) -- (4,0);

  \draw[black!50, semithick, densely dashed] (1,0) -- (2,1) -- (3,0);

  \fill[black!10] (0,0) -- (2,2) -- (4,0) -- cycle;

  \fill[black!22] (1,0) -- (2,1) -- (3,0) -- cycle;

  \draw[black, semithick] (0,0) -- (2,2) -- (4,0);
  \draw[black!50, semithick, densely dashed] (1,0) -- (2,1) -- (3,0);

  \node[ptlabel, above] at (2,2) {$\lambda_1$};
  \node[ptlabel, above] at (2,1.08) {$\lambda_2$};

  \node[ptlabel, left=2pt] at (0.9, 0.90) {$\Lambda_1$};
  \node[ptlabel, right=2pt] at (2.6, 0.55) {$\Lambda_2$};
\end{scope}
\node[paneltitle] at (10.35, -0.60) {(b) Tent functions and landscape layers};

\end{tikzpicture}%
}

%% file: experiments/results/table_datasets.tex
\begin{tabular}{lrrrcl}
\toprule
Dataset & $n$ & $T$ & \#cls & Imbalance & Domain \\
\midrule
Coffee & 56 & 286 & 2 & 1.07 & Spectroscopy \\
GunPoint & 200 & 150 & 2 & 1.00 & Motion \\
ECGFiveDays & 884 & 136 & 2 & 1.00 & ECG \\
PowerCons & 360 & 144 & 2 & 1.00 & Energy \\
ToeSegmentation1 & 268 & 277 & 2 & 1.09 & Gait \\
ECG200 & 200 & 96 & 2 & 1.99 & ECG \\
FordA & 4921 & 500 & 2 & 1.06 & Engine sensor \\
FordB & 4446 & 500 & 2 & 1.03 & Engine sensor \\
ArrowHead & 211 & 251 & 3 & 1.25 & Shape outline \\
Meat & 120 & 448 & 3 & 1.00 & Spectroscopy \\
Plane & 210 & 144 & 7 & 1.00 & Shape outline \\
Trace & 200 & 275 & 4 & 1.00 & Sensor transient \\
\bottomrule
\end{tabular}

%% file: experiments/results/table_params.tex
\begin{tabular}{lrrrr}
\toprule
Dataset & \multicolumn{2}{c}{CGSSN/$k$-NN} & \multicolumn{2}{c}{OPN} \\
\cmidrule(lr){2-3} \cmidrule(lr){4-5}
 & $\tau$ & $n$ & $\tau$ & $n$ \\
\midrule
Coffee & 13 & 2 & 5 & 5 \\
GunPoint & 7 & 3 & 2 & 4 \\
ECGFiveDays & 6 & 3 & 4 & 4 \\
PowerCons & 5 & 4 & 3 & 5 \\
ToeSegmentation1 & 10 & 3 & 2 & 5 \\
ECG200 & 6 & 2 & 2 & 4 \\
FordA & 6 & 3 & 3 & 5 \\
FordB & 5 & 3 & 5 & 5 \\
ArrowHead & 11 & 2 & 5 & 3 \\
Meat & 12 & 3 & 5 & 6 \\
Plane & 8 & 3 & 4 & 5 \\
Trace & 3 & 3 & 2 & 5 \\
\bottomrule
\end{tabular}

%% file: experiments/results/table_extraction_time.tex
\begin{tabular}{lrrrrr}
\toprule
Dataset & HVG & NVG & OPN & CGSSN & KNN \\
\midrule
Coffee & 3.8 & 2.8 & 1.2 & 0.1 & 0.8 \\
GunPoint & 0.2 & 0.2 & 0.1 & 0.1 & 0.3 \\
ECGFiveDays & 0.6 & 0.6 & 0.2 & 0.2 & 0.6 \\
PowerCons & 0.3 & 0.3 & 0.2 & 4.7 & 0.3 \\
ToeSegmentation1 & 0.8 & 1.3 & 0.2 & 0.3 & 1.3 \\
ECG200 & 0.1 & 0.1 & 0.1 & 0.1 & 0.1 \\
FordA & 33.5 & 46.3 & 2.1 & 10.5 & 52.8 \\
FordB & 29.7 & 42.2 & 3.6 & 10.1 & 48.4 \\
ArrowHead & 0.5 & 0.8 & 0.1 & 0.1 & 1.1 \\
Meat & 1.5 & 2.2 & 0.2 & 0.1 & 4.7 \\
Plane & 0.2 & 0.2 & 0.1 & 0.1 & 0.2 \\
Trace & 0.8 & 0.9 & 0.3 & 0.1 & 1.3 \\
\midrule
\textbf{Total} & \textbf{71.8} & \textbf{97.9} & \textbf{8.4} & \textbf{26.4} & \textbf{111.8} \\
\bottomrule
\end{tabular}

%% file: experiments/results/table_f1.tex
\begin{tabular}{lrrrrr}
\toprule
Dataset & HVG & NVG & OPN & CGSSN & KNN \\
\midrule
ArrowHead & 0.487\scriptsize{$\pm$0.045} & 0.185\scriptsize{$\pm$0.002} & 0.472\scriptsize{$\pm$0.057} & 0.664\scriptsize{$\pm$0.069} & \textbf{0.692}\scriptsize{$\pm$0.058} \\
Coffee & 0.341\scriptsize{$\pm$0.018} & 0.341\scriptsize{$\pm$0.018} & \textbf{0.928}\scriptsize{$\pm$0.040} & 0.709\scriptsize{$\pm$0.176} & 0.629\scriptsize{$\pm$0.059} \\
ECG200 & 0.457\scriptsize{$\pm$0.055} & 0.399\scriptsize{$\pm$0.005} & 0.492\scriptsize{$\pm$0.034} & \textbf{0.618}\scriptsize{$\pm$0.084} & 0.504\scriptsize{$\pm$0.064} \\
ECGFiveDays & 0.627\scriptsize{$\pm$0.032} & 0.332\scriptsize{$\pm$0.001} & 0.686\scriptsize{$\pm$0.036} & \textbf{0.846}\scriptsize{$\pm$0.023} & 0.568\scriptsize{$\pm$0.016} \\
FordA & \textbf{0.992}\scriptsize{$\pm$0.003} & 0.339\scriptsize{$\pm$0.000} & 0.754\scriptsize{$\pm$0.010} & 0.710\scriptsize{$\pm$0.016} & 0.655\scriptsize{$\pm$0.025} \\
FordB & 0.337\scriptsize{$\pm$0.000} & 0.337\scriptsize{$\pm$0.000} & \textbf{0.770}\scriptsize{$\pm$0.007} & 0.685\scriptsize{$\pm$0.010} & 0.657\scriptsize{$\pm$0.009} \\
GunPoint & 0.584\scriptsize{$\pm$0.090} & 0.333\scriptsize{$\pm$0.000} & 0.654\scriptsize{$\pm$0.060} & 0.659\scriptsize{$\pm$0.038} & \textbf{0.769}\scriptsize{$\pm$0.048} \\
Meat & 0.286\scriptsize{$\pm$0.085} & 0.167\scriptsize{$\pm$0.000} & \textbf{0.737}\scriptsize{$\pm$0.068} & 0.714\scriptsize{$\pm$0.032} & 0.584\scriptsize{$\pm$0.049} \\
Plane & 0.036\scriptsize{$\pm$0.000} & 0.036\scriptsize{$\pm$0.000} & 0.736\scriptsize{$\pm$0.041} & \textbf{0.845}\scriptsize{$\pm$0.042} & 0.843\scriptsize{$\pm$0.076} \\
PowerCons & 0.651\scriptsize{$\pm$0.070} & 0.333\scriptsize{$\pm$0.000} & 0.825\scriptsize{$\pm$0.088} & \textbf{0.875}\scriptsize{$\pm$0.040} & 0.779\scriptsize{$\pm$0.075} \\
ToeSegmentation1 & 0.342\scriptsize{$\pm$0.005} & 0.343\scriptsize{$\pm$0.002} & \textbf{0.739}\scriptsize{$\pm$0.051} & 0.650\scriptsize{$\pm$0.034} & 0.622\scriptsize{$\pm$0.042} \\
Trace & 0.100\scriptsize{$\pm$0.000} & 0.100\scriptsize{$\pm$0.000} & 0.569\scriptsize{$\pm$0.077} & \textbf{0.979}\scriptsize{$\pm$0.028} & 0.681\scriptsize{$\pm$0.083} \\
\midrule
\#wins & 1 & 0 & 4 & 5 & 2 \\
\bottomrule
\end{tabular}

%% file: experiments/results/table_auc.tex
\begin{tabular}{lrrrrr}
\toprule
Dataset & HVG & NVG & OPN & CGSSN & KNN \\
\midrule
ArrowHead & 0.684\scriptsize{$\pm$0.066} & 0.500\scriptsize{$\pm$0.000} & 0.652\scriptsize{$\pm$0.101} & 0.812\scriptsize{$\pm$0.049} & \textbf{0.853}\scriptsize{$\pm$0.017} \\
Coffee & 0.500\scriptsize{$\pm$0.000} & 0.500\scriptsize{$\pm$0.000} & \textbf{0.994}\scriptsize{$\pm$0.012} & 0.764\scriptsize{$\pm$0.128} & 0.711\scriptsize{$\pm$0.103} \\
ECG200 & 0.529\scriptsize{$\pm$0.030} & 0.500\scriptsize{$\pm$0.000} & 0.527\scriptsize{$\pm$0.056} & \textbf{0.694}\scriptsize{$\pm$0.113} & 0.495\scriptsize{$\pm$0.099} \\
ECGFiveDays & 0.636\scriptsize{$\pm$0.040} & 0.500\scriptsize{$\pm$0.000} & 0.749\scriptsize{$\pm$0.032} & \textbf{0.932}\scriptsize{$\pm$0.020} & 0.618\scriptsize{$\pm$0.034} \\
FordA & \textbf{0.994}\scriptsize{$\pm$0.003} & 0.500\scriptsize{$\pm$0.000} & 0.833\scriptsize{$\pm$0.008} & 0.777\scriptsize{$\pm$0.008} & 0.705\scriptsize{$\pm$0.025} \\
FordB & 0.498\scriptsize{$\pm$0.002} & 0.500\scriptsize{$\pm$0.000} & \textbf{0.854}\scriptsize{$\pm$0.004} & 0.750\scriptsize{$\pm$0.013} & 0.714\scriptsize{$\pm$0.011} \\
GunPoint & 0.590\scriptsize{$\pm$0.086} & 0.500\scriptsize{$\pm$0.000} & 0.694\scriptsize{$\pm$0.082} & 0.730\scriptsize{$\pm$0.081} & \textbf{0.852}\scriptsize{$\pm$0.043} \\
Meat & 0.558\scriptsize{$\pm$0.031} & 0.500\scriptsize{$\pm$0.000} & \textbf{0.896}\scriptsize{$\pm$0.038} & 0.886\scriptsize{$\pm$0.046} & 0.770\scriptsize{$\pm$0.049} \\
Plane & 0.500\scriptsize{$\pm$0.000} & 0.500\scriptsize{$\pm$0.000} & 0.961\scriptsize{$\pm$0.011} & \textbf{0.980}\scriptsize{$\pm$0.007} & 0.980\scriptsize{$\pm$0.011} \\
PowerCons & 0.744\scriptsize{$\pm$0.075} & 0.500\scriptsize{$\pm$0.000} & 0.914\scriptsize{$\pm$0.052} & \textbf{0.952}\scriptsize{$\pm$0.032} & 0.837\scriptsize{$\pm$0.072} \\
ToeSegmentation1 & 0.493\scriptsize{$\pm$0.011} & 0.500\scriptsize{$\pm$0.000} & \textbf{0.831}\scriptsize{$\pm$0.058} & 0.728\scriptsize{$\pm$0.064} & 0.675\scriptsize{$\pm$0.059} \\
Trace & 0.500\scriptsize{$\pm$0.000} & 0.500\scriptsize{$\pm$0.000} & 0.876\scriptsize{$\pm$0.029} & \textbf{0.999}\scriptsize{$\pm$0.001} & 0.900\scriptsize{$\pm$0.029} \\
\midrule
\#wins & 1 & 0 & 4 & 5 & 2 \\
\bottomrule
\end{tabular}